\documentclass[11pt]{article}
\usepackage{amstext,amssymb,amsmath,amsbsy}

\usepackage{hyperref}
\usepackage{amscd}
\usepackage{amsfonts}
\usepackage{indentfirst}
\usepackage{verbatim}
\usepackage{amsmath}
\usepackage{amsthm}
\usepackage{enumerate}
\usepackage{graphicx}
\usepackage{subfig}
\usepackage{color}
\usepackage[OT1]{fontenc}
\usepackage[latin1]{inputenc}
\usepackage[english]{babel}
\usepackage{amssymb}

\newtheorem{theorem}{Theorem}
\newtheorem{lemma}{Lemma}

\newtheorem{proposition}{Proposition}
\newtheorem{corollary}{Corollary}
\newtheorem{remark}{Remark}

\newtheorem{definition}{Definition}
\numberwithin{equation}{section}

\newcommand{\proofend}{\hfill $\Box$ }

\newcommand{\sge}{\gtrsim}

\newcommand{\supp}{\operatorname{supp}}
\newcommand{\dist}{\operatorname{dist}}

\newcommand{\dive}{\operatorname{div}}

\newcommand{\ff}{ {\bf \varphi}}

\newcommand{\tr}{ {\mbox{trace}}}

\newcommand{\eps}{\varepsilon}

\newcommand{\loc}{_{loc}}

\newcommand{\mN}{\mathbb{N}}

\newcommand{\mR}{\mathbb{R}}

\newcommand{\bw}{{\bf w}}

\newcommand{\bH}{{\bf H}}

\title{Limiting absorption principle and well-posedness for the Helmholtz equation with sign changing coefficients}

\author{Hoai-Minh Nguyen \footnote{EPFL SB MATHAA CAMA, Station 8,  CH-1015 Lausanne, hoai-minh.nguyen@epfl.ch} 
}

\begin{document}

\maketitle 

\begin{abstract}
In this paper,  we  investigate the limiting absorption principle associated to and  the well-posedness of the Helmholtz equations with sign changing coefficients which are used to model  negative index materials. Using the reflecting technique introduced in \cite{Ng-Complementary}, we first derive Cauchy problems from these equations. The limiting absorption principle and the well-posedness are then obtained via various a priori estimates for these Cauchy problems. There approaches are proposed to obtain the a priori estimates. The first one  follows from a priori estimates of elliptic  systems equipped general complementing boundary conditions due to  Agmon, Douglis, and Nirenberg in their classic work \cite{ADNII}. The second approach, which is  complement to the first one, is variational and based on the Dirichlet principle. The last approach, which is complement to the second one,  is also variational and uses the multiplier technique. Using these approaches, we are able to obtain new results on the well-posedness of these equations for which the conditions on the coefficients are imposed partially or not strictly on the interface of sign changing coefficients. In particular, the well-posedness can hold even in the case the contrast of the coefficients across the sign changing interfaces is arbitrary. This allows us to rediscover and extend  known results obtained by the integral method, the pseudo differential operator theory, and the T-coercivity approach.  The unique  solution, obtained by the limiting absorption principle,  is {\bf not}  in $H^1_{\loc}(\mR^d)$ as usual and possibly  {\bf not even} in  $L^2_{\loc}(\mR^d)$. The optimality of our results is  also discussed. 

\end{abstract}



\tableofcontents

\section{Introduction}

This paper deals with  the Helmholtz equation with sign changing coefficients which are used to model  negative index materials (NIMs). NIMs were first investigated theoretically by Veselago in \cite{Veselago}. The  existence of such materials was confirmed by Shelby,  Smith, and Schultz in \cite{ShelbySmithSchultz}. The study of NIMs has attracted a lot attention in the scientific community thanks to their many possible applications such as superlensing and cloaking using complementary media, and cloaking  a source via anomalous localized resonance. 

\medskip
We next mention briefly these three applications of NIMs. Superlensing using negative index materials  was suggested by Veselago in \cite{Veselago} for a slab lens (a slab of index $-1$) using the ray theory. 
Later, cylindrical lenses in the two dimensional  quasistatic regime, the Veselago slab lens  and cylindrical lenses in the finite frequency regime, and   spherical lenses in the finite frequency regime were studied by Nicorovici, McPhedran, and Milton  in \cite{NicoroviciMcPhedranMilton94},  Pendry in \cite{PendryNegative, PendryCylindricalLenses}, and  Pendry and Ramakrishna in \cite{PendryRamakrishna} respectively for constant  isotropic objects. Superlensing using NIMs (or more precisely using complementary media) for arbitrary objects in the acoustic and electromagnetic settings  was established in  \cite{Ng-Superlensing, Ng-Superlensing-Maxwell} for schemes inspired by \cite{NicoroviciMcPhedranMilton94, PendryCylindricalLenses, PendryRamakrishna}  and guided by the concept of reflecting complementary media introduced and studied in \cite{Ng-Complementary}.  Cloaking using  complementary media   was suggested and investigated numerically by Lai et al. in \cite{LaiChenZhangChanComplementary}. Cloaking an arbitrary inhomogeneous object using complementary media  was proved in \cite{Ng-Negative-cloaking} for  the  quasi-static regime and later extended in \cite{MinhLoc2} for the finite frequency regime. The schemes used there are inspired by  \cite{LaiChenZhangChanComplementary} and \cite{Ng-Complementary}.  
Cloaking  a source via anomalous localized resonance was discovered by Milton and Nicorovici for constant symmetric plasmonic structures in the two dimensional
quasistatic regime in \cite{MiltonNicorovici} (see also \cite{MiltonNicoroviciMcPhedranPodolskiy, NicoroviciMcPhedranMilton94}) for dipoles.  Cloaking an arbitrary source concentrated on a manifold of codimension 1  in an arbitrary medium via anomalous localized resonance was proposed and established in \cite{Ng-CALR-CRAS, Ng-CALR, Ng-CALR-finite}.  Other contributions  are  \cite{AmmariCiraoloKangLeeMilton, AmmariCiraoloKangLeeMilton2, BouchitteSchweizer10, KohnLu, MinhLoc1}  in which special structures and partial aspects were investigated. A survey on the mathematics progress of these applications can be  found in \cite{Ng-Negative-Review}. 
It is worthy to note that in the applications of NIMs mentioned above, the localized resonance,  i.e., the field blows up in some regions and remains bounded in some others as the loss goes to 0, might appear.  

\medskip

In this paper, we investigate the well-posedness of the Helmholtz equation with sign changing coefficients: the stability aspect. To ensure to obtain 
 physics solutions, we also study  the limiting absorption principle associated to this equation.
 Let $k > 0$ and let $A$ be a (real) uniformly elliptic symmetric matrix defined on $\mR^d$ ($d \ge 2$), and $\Sigma$ be a  bounded real function defined on $\mR^d$.  Assume that  
\begin{equation*}
A(x) = I \mbox{ in } \mR^d \setminus B_{R_0},  \quad A \mbox{ is piecewise } C^1, 
\end{equation*}
and 
\begin{equation*}
\Sigma(x) = 1 \mbox{ in } \mR^d \setminus B_{R_0}, 
\end{equation*}
for  some  $R_0 > 0$. 
Here and in what follows, for $R>0$,  $B_R$ denotes the ball in $\mR^d$ centered at the origin and  of radius $R$.  Let $D \subset \subset B_{R_0}$ be a bounded open subset in $\mR^d$  of class $C^2$. 
Set, for $\delta \ge 0$,  
\begin{equation}\label{def-sd}
s_\delta (x) = \left\{\begin{array}{cl} - 1 - i \delta & \mbox{ in } D, \\[6pt]
1 & \mbox{ in } \mR^d \setminus D. 
\end{array}\right.
\end{equation}
We are interested in the well-posedness in the class of outgoing solutions of the following equation 
\begin{equation}\label{Main-eq}
\dive(s_0A \nabla u_0) + k^2 s_0 \Sigma u_0 = f \mbox{ in } \mR^d, 
\end{equation} 
and the limiting absorption principle associated with it, i.e.,  the  convergence of   $u_\delta$ to $u_0$ (in an appropriate sense) under various conditions on $A$ and $\Sigma$. Here $u_\delta \in H^1(\mR^d)$ $(\delta > 0)$ is the unique solution of  the equation
\begin{equation}\label{Main-eq-delta}
\dive(s_\delta A \nabla u_\delta) + k^2 s_0 \Sigma u_\delta + i \delta u_\delta  = f \mbox{ in } \mR^d. 
\end{equation} 
Recall that a solution $v \in H^1_{\loc}(\mR^d \setminus B_{R})$ of the equation
\begin{equation*}
\Delta v + k^2 v = 0 \mbox{ in } \mR^d \setminus B_{R}, 
\end{equation*} 
for some $R>0$, is said to satisfy the outgoing condition if 
\begin{equation*}
\partial_r v - i k v = o(r^{-\frac{d-1}{2}}) \mbox{ as } r = |x| \to + \infty. 
\end{equation*}
Physically, $k$ is the frequency, $(s_\delta A, s_0 \Sigma)$ is the material parameter of the medium, and $\delta$ describes the loss of the material.  We denote
$$
\Gamma = \partial D, 
$$
and, for $\tau > 0$, 
\begin{equation}\label{Otau}
D_\tau = \big\{x \in D;\;  \dist(x, \Gamma) < \tau \big\}
\end{equation}
\begin{equation}\label{O-tau}
D_{-\tau} = \big\{x \in \mR^d \setminus \bar D; \; \dist(x, \Gamma)  < \tau \big\}. 
\end{equation}
As usual, $\bar D$ denotes the closure of $D$ for a subset $D$ in $\mR^d$.  

\medskip
The well-posedness of the Helmholtz equation
with sign changing coefficients was first established by  Costabel and Stephan in \cite{CostabelErnst}. They proved, by the integral approach, that \eqref{Main-eq} is well-posed if $A = I$ in $\mR^d \setminus D$ and $A = \lambda I$ in $D$ provided that  $\lambda$ is a positive constant not equal to 1.  Later, Ola in \cite{Ola} proved, using the integral method and the pseudo-diiferential operators theory, that  \eqref{Main-eq} is well-posed in three and higher dimensions if $\Gamma$ is strictly convex and  connected even though $\lambda =1$, i.e.,  $A= I$ in $\mR^d$. His result was  extended for the case where $\Gamma$ has two strictly convex connected components by Kettunen, Lassas,  and Ola in \cite{KettunenLassas}.  Recently, the well-posedness  was extensively studied by Bonnet-Ben Dhia, Ciarlet, and their coauthors in \cite{AnneSophieChesnelCiarlet1, AnneSophieChesnelCiarlet3, AnneSophieChesnelCiarlet2, AnneSophieChesnelCiarlet1-1, AnneSophie-Ciarlet0, AnneSophie-Ciarlet01, ChesnelCiarlet} by  T-coercivity approach. This approach was introduced by Bonnet-Ben Dhia,  Ciarlet, and Zw\"olf  in \cite{AnneSophie-Ciarlet0} and  is related to the (Banach-Necas-Babuska) inf-sup
condition. The sharpest result for the acoustic setting in this direction, obtained by Bonnet-Ben Dhia, Chesnel, and Ciarlet in \cite{AnneSophieChesnelCiarlet1},  is  that   \eqref{Main-eq} is well-posed in the Fredholm sense in $H^1$ (this means that the compactness holds \footnote{They considered the bounded setting and the uniqueness is not ensured in general.}), if $A$ is  isotropic, i.e., $A =  a I$ for some positive function $a$,  and  the contrast of $a$ is not 1 on each connected component of $\Gamma$. 

\medskip
In this paper, we are interested in  the limiting absorbtion principle and the  well-posedness of  \eqref{Main-eq} for solutions obtained by the limiting absorption process. Our starting point  is to obtain Cauchy's problems using the reflecting technique introduced in \cite{Ng-Complementary}. The idea is simple as follows. Let $F: U \setminus \bar D \to D_\tau$ be a reflection through $\Gamma$, i.e., $F$ is a diffeomorphism and $F(x) = x$ on $\Gamma $ for some smooth open set $D \subset \subset U$ and for some $\tau > 0$. Set $v_\delta = u_\delta \circ F^{-1}$. By a change of variables (see Lemma~\ref{lem-TO}), it follows from \eqref{def-sd} that  
\begin{equation*}
\dive(F_*A \nabla v_\delta)  + k^2 F_*\Sigma v_\delta = F_* f  + O(\delta v_\delta) \mbox{ in } D_\tau, 
\end{equation*} 
\begin{equation*}
\dive( \nabla A u_\delta)  + k^2 \Sigma u_\delta = s_0^{-1}f + O(\delta u_\delta) + O(\delta f) \mbox{ in } D_\tau, 
\end{equation*}
\begin{equation*}
v_\delta - u_\delta = 0 \mbox{ on } \Gamma \quad \mbox{ and } F_*A \nabla v_\delta \cdot \nu - A \nabla u_\delta \cdot \nu = i \delta A \nabla u_\delta \cdot \nu \mbox{ on } \Gamma. 
\end{equation*} 
Here and in what follows, for a matrix $a$, a function $\sigma$, and a diffeomorphism $T$, the following standard notations are used: 
\begin{equation}\label{TO}
T_*a(y) = \frac{D T  (x)  a(x) D T ^{T}(x)}{J(x)} \quad  \mbox{  and  }  \quad  T_*\sigma(y) = \frac{\sigma(x)}{J(x)}, 
\end{equation}
where 
\begin{equation*}
 J(x) = |\det D T(x)| \quad  \mbox{ and } \quad x =T^{-1}(y). 
\end{equation*}
Here we  denote $O(v)$ a quantity whose $L^2$-norm is bounded by $C \| v\|_{L^2}$ for some positive constant $C$ independent of $\delta$  and $v$ for $0  < \delta < 1$. 
We hence obtain Cauchy's problems  for $(u_\delta, v_\delta)$ in $D_\tau$  by considering $O(\delta v_\delta)$, $O(\delta u_\delta)$, $O(\delta f)$, and $i \delta A \nabla u_\delta \cdot \nu$ like  given data which are formally 0 if $\delta = 0$. The use of reflections to study NIMs was also considered by Milton et al. in \cite{Milton-folded} and by  Bonnet-Ben Dhia, Ciarlet, and their coauthors in their T-coercivity approach. However,
there is a difference between the use of reflections in \cite{Milton-folded},  in the  T-coercivity approach,  and in our work. In \cite{Milton-folded}, the
authors used reflections as a change of variables to obtain a new simple setting from an old
more complicated one and hence the analysis of the old problem becomes simpler. In the  T-coercivity approach, the authors used a standard reflection to construct test functions for  the inf-sup condition to obtain an a priori estimate for the solution. Our use of  reflections is to derive the Cauchy problems. This can be done in a very flexible way via a change of variables formula stated in Lemma~\ref{lem-TO} as observed in \cite{Ng-Complementary}. The limiting absorption principle and the  well-posedness of \eqref{Main-eq} are then based on a priori estimates for these Cauchy problems under various conditions on $A, \Sigma$, $F_*A$, and $F_*\Sigma$ in $D_\tau$.  Appropriate choices of reflections are important in the applications and  discussed later (Corollaries~\ref{cor1}, \ref{cor2}, and \ref{cor3}). 

\medskip
In this paper, we introduce three  approaches to obtain  a priori estimates for the Cauchy problems.  The first one follows  from a priori estimates for {\bf elliptic  systems} imposing general {\bf implementing boundary conditions} (see Definition~\ref{Def1}) due to Agmon, Douglis, and Nirenberg in  their classic work \cite{ADNII}. Applying their result, we can prove in Section~\ref{sect-Thm0}: 
\begin{itemize}
\item[1.] Assume that  $A_+ := A \big|_{\mR^d \setminus D} \in C^1(\bar D_{-\tau})$ and $A_-: = A \big|_{D} \in C^1(\bar D_\tau)$ for some small positive constant $\tau$,  and  $A_{+}$ and $A_{-}$ satisfy the (Cauchy) complementing boundary condition  on $\Gamma$. Then the limiting absorption principle  and the well-posedness  in $H^1_{\loc}(\mR^d)$ for \eqref{Main-eq} hold (Theorem~\ref{thm0} in Section~\ref{sect-Thm0}).  
\end{itemize}
In fact, we  establish that the  conclusions  hold if $F_*A_+$ and $A_-$ satisfy the (Cauchy) complementing boundary  condition  on $\Gamma$ where  $F$ is the standard reflection in  \eqref{def-reflection}.  Using the characterization of complementing boundary condition established  in Proposition~\ref{pro-complementing}, one can prove that $F_*A_+$ and $A_-$ satisfy the (Cauchy) complementing boundary condition  on $\Gamma$ if and only if $A_+$ and $A_-$ do;  this implies the first result.  Using the first result, one obtains new conditions for which  the well-posedness and  the limiting absorption principle hold. In particular, the condition $A_+ > A_-$ or $A_- > A_+$ on each connected component of $\Gamma$ is sufficient for the conclusion (see Corollary~\ref{cor0}). 
Here and in what follows, we use the following standard notation for a matrix $M$:  $M > 0$ means that $\langle M x , x \rangle  > 0 $ for all $x \neq 0$ where $\langle \cdot, \cdot \rangle$ denotes the Euclidean scalar product in $\mR^d$.  To our knowledge, Corollary~\ref{cor0} is new and cannot be obtained using the known approaches mentioned above.  Corollary~\ref{cor0} is in the same spirit of the one of Bonnet-Ben Dhia, Chesnel, and Ciarlet in \cite{AnneSophieChesnelCiarlet1}; nevertheless, $A_+$ and $A_-$ are not assumed to be isotropic here.
One can verify that if $F_*A_+  = A_-$ on $\Gamma$ then the complementing boundary condition is not satisfied (see Proposition~\ref{pro-complementing}). To deal with this situation, we develop a second approach to obtain a priori estimates for the Cauchy problems in Section~\ref{sect-Thm1}. This approach is variational  and based on the  Dirichlet principle. 
Using this approach, we can establish:
\begin{itemize}
\item[2.] Assume that there exist $\tau > 0$ (small), a smooth open set $U \supset \supset D$, and 
a reflection $F: U \setminus D \to D_\tau$, i.e., $F$ is diffeomorphism and $F(x) = x$ on $\Gamma$,  such that, on every connected component of $D_\tau$, 
\begin{equation}\label{T1-cond*}
\mbox{either }  \quad A - F_*A \sge  \dist(x, \Gamma)^{\alpha} I \quad  \mbox{ or } \quad 
F_*A - A \sge \dist(x, \Gamma)^{\alpha} I,  
\end{equation}
 for some $0 \le \alpha < 2$. Then the limiting absorption principle  and the well-posedness for \eqref{Main-eq} hold (Theorem~\ref{thm1} in Section~\ref{sect-Thm1}). 
\end{itemize}

The  unique solution, which is obtained by the limiting absorption principle,  might  {\bf not} be in $H^1_{\loc}(\mR^d)$ in this case;  the  proof of the uniqueness is nonstandard. 
The appropriate space in which the solution is defined  is revealed by the limiting absorption principle; more precisely, by a priori estimates obtained for $u_\delta$ defined in \eqref{Main-eq-delta}.  Once the uniqueness is obtained, the stability is based on a compactness argument. A new compactness criterion in $L^2$ (Lemma~\ref{lem3}) is established in this process and the condition $\alpha < 2$ is required there. Various consequences of this result are given in Section~\ref{sect-Thm1} (Corollaries~\ref{cor1} and \ref{cor2}). The choice of the reflections is crucial in deriving these consequences).  Theorem~\ref{thm1} implies, unifies,  and extends the known results mentioned above.
In particular, using Corollary~\ref{cor2}, one can derive that  the well-posedness holds under the condition that  $\Gamma$ is strictly convex (the number of connected component of $\Gamma$ is not imposed) and $A$ is  isotropic and constant {\bf only}  on each small connected component neighborhood of $\Gamma$ in three and higher dimensions.  A variant of the result of Ola in \cite{Ola} in two dimensions holds and is also contained in Theorem~\ref{thm1}. 


\medskip
Similar conclusion still holds in the  case $F_*A = A$ in $D_\tau$  under additional assumptions on $\Sigma$ and $F_*\Sigma$ in $D_\tau$. To reach the conclusion in this case, we propose a third approach to deal with the Cauchy problems in Section~\ref{sect-Thm2}. It is   variational  and based on the multiplier technique.   In this direction, we can prove the following result:
\begin{itemize}
\item[3.] Assume that there exist $\tau > 0$ (small), a smooth open subset $U \supset \supset D$, and 
a reflection $F: U \setminus D \to D_\tau$ such that either 
\begin{equation}\label{T2-cond*-1}
F_*A - A \ge 0 \quad  \mbox{ and } \quad \Sigma - F_*\Sigma  \sge  \dist(x, \Gamma)^\beta 
\end{equation}
or 
\begin{equation}\label{T2-cond*-2}
A - F_*A \ge 0 \quad   \mbox{ and } \quad  F_*\Sigma - \Sigma   \sge  \dist(x, \Gamma)^\beta,    
\end{equation}
in each connected component of $D_\tau$ for some $\beta > 0$. Then the limiting absorption principle  and the well-posedness for \eqref{Main-eq} hold (Theorem~\ref{thm2} in Section~\ref{sect-Thm2}).  
\end{itemize}
The unique solution $u$, in this case,  is {\bf not even} in  $L^2_{\loc}(\mR^d)$ and $f$ is assumed  to be 0 near $\Gamma$. The appropriate space for which the solution is defined is again revealed  by the limiting absorption principle. Once the uniqueness is established, the stability is based on a compactness argument. Due to the lack of $L^2$-control, the  compactness argument used in this case is non-standard and different from the one used in the second setting (see the proofs of Theorem~\ref{thm1} and Theorem~\ref{thm2}). A simple application of this result is given in Corollary~\ref{cor3} which is a supplement to 
Corollary~\ref{cor2} in two dimensions. As far as we know, Theorem~\ref{thm2} is the first result on the limiting absorption principle and the well-posedness for the Helmholtz equations with sign changing coefficients where the conditions on the coefficients contain the zero order term $\Sigma$.

\medskip
It is known that in the case $(F_*A, F_*\Sigma) = (A, \Sigma)$ in $D_\tau$, the  localized resonance might appear. Media with this property are roughly speaking called reflecting complementary media introduced and studied in \cite{Ng-Complementary, Ng-Superlensing-Maxwell} for the Helmholtz and Maxwell equations respectively. The notion of reflecting complementary media plays an important roles in various applications of NIMs  mentioned previously  as discussed in  \cite{Ng-Superlensing, Ng-CALR-CRAS, Ng-CALR, Ng-Negative-cloaking, Ng-Superlensing-Maxwell, MinhLoc1, MinhLoc2}.  The results obtained in this paper, in particular from the second and the third results,  showed that the complementary property of media is necessary for the occurrence of the resonance.  In Section~\ref{sect-resonant}, we show that even in the case  $(F_*A, F_*\Sigma) = (A, \Sigma)$ in $B(x_0, r_0) \cap D_\tau$ for some  $x_0 \in \Gamma$ and $r_0 > 0$,  the system is resonant in the following sense (see Proposition~\ref{pro-resonant}): There exists $f$ with $\supp f \subset \subset B_{R_0} \setminus \Gamma$  such that  $\limsup_{\delta \to 0}\| u_{\delta}\|_{L^2(K)} = + \infty$ for some  $K \subset \subset B_{R_0} \setminus  \Gamma$.  Here and in what follows $B(x, r)$ denotes the ball centered at $x$ and of radius $r$. 
This also implies the optimality of the results mentioned above. The proof of Proposition~\ref{pro-resonant} is based on a three sphere inequality and has roots from \cite{Ng-CALR}.



%


\medskip
The paper is organized as follows. Sections~\ref{sect-Thm0}, \ref{sect-Thm1}, and \ref{sect-Thm2} are devoted to the proof of the three main results mentioned above and their consequences respectively. In Section~\ref{sect-resonant}, we disscuss the optimality of these results.  


\section{An approach via  a priori estimates  of elliptic systems imposed complementing boundary conditions}\label{sect-Thm0}

A useful  simple technique suggested to study  the Helmholtz equations with sign changing coefficients is the reflecting one introduced in \cite{Ng-Complementary}.  Applying this  technique,  we obtain Cauchy problems from the Helmholtz equations with sign changing coefficients. An important part in the investigation of the  well-posedness and the limiting absorption principle   is then to obtain appropriate a priori estimates for these Cauchy problems. In this section, these follow from an  estimate near the boundary of solutions of {\bf elliptic  systems} imposed Cauchy data due to Agmon, Douglis, and Nirenberg in their classic work \cite{ADNII} (see also \cite{Lo}).  Before stating the result, let us recall the notation of complementing boundary condition with respect to the Cauchy data derived from \cite{ADNII}. 

\begin{definition}[Agmon, Douglis, Nirenberg \cite{ADNII}]\label{Def1} Two constant positive symmetric matrices $A_1$ and $A_2$ are said to satisfy the (Cauchy) complementing boundary condition with respect to  direction $e  \in \partial B_1$ if and only if for all $\xi \in \mR^d_{e, 0} \setminus \{0 \}$, the only solution $(u_1(x), u_2(x))$ of the form $\big(e^{i \langle y,  \xi \rangle  } v_1(t), e^{i \langle y,  \xi \rangle} v_2(t) \big)$ with $x  =  y + t e$ where  $t = \langle x, e \rangle $, of the  following system 
\begin{equation*}
\left\{ \begin{array}{c}
\dive(A_1 \nabla u_1) = \dive(A_2 \nabla u_2) = 0 \mbox{ in } \mR^d_{e, +}, \\[6pt]
u_1 = u_2 \mbox{ and } A_1 \nabla u_1 \cdot e = A_2 \nabla u_2 \cdot e \mbox{ on } \mR^d_{e, 0}, 
\end{array}\right.
\end{equation*}
which is bounded in $\mR^d_{e, +}$ is $(0, 0)$. 
\end{definition}

Here and in what follows, for a unit vector  $e \in \mR^d$, the following notations are used
\begin{equation}\label{notation-Rd}
\mR^d_{e, +} = \{\xi \in \mR^d; \; \langle \xi, e \rangle > 0 \} \quad \mbox{ and } \quad \mR^d_{e, 0} = \{\xi \in \mR^d; \; \langle \xi, e \rangle = 0 \}. 
\end{equation}
Recall that  $\langle \cdot, \cdot \rangle$ denotes the Euclidean scalar product in $\mR^d$. 

\medskip

We are ready to state the main result of this section: 

\begin{theorem}\label{thm0}  Let $f \in L^2(\mR^d)$ with $\supp f \subset \subset B_{R_0}$, and  
let $u_\delta \in H^1(\mR^d)$ $(0 < \delta < 1)$ be the unique solution of \eqref{Main-eq-delta}. Assume that $A_+: = A  \big|_{\mR^d \setminus \bar D} \in C^1(\bar D_{-h})$ and  $A_-: = A \big|_{D} \in C^1(\bar D_h )$, and $A_+(x), A_-(x)$ satisfy the (Cauchy) complementing boundary condition with respect to direction  $\nu(x)$  for all $x \in \Gamma$.  Then 
\begin{equation}\label{T0-conclusion0}
\| u_\delta\|_{H^1(B_R)} \le C_R \| f\|_{L^2(\mR^d)} \quad \forall \, R > 0, 
\end{equation}
for some positive constant $C_R$ independent of $\delta$ and $f$.  Moreover,  $u_\delta \to u_0$ weakly in $H^1_{\loc}(\mR^d)$, as $\delta \to 0$, where  $u_0 \in H^1_{\loc}(\mR^d)$ is the  unique outgoing solution of \eqref{Main-eq}.  We also have 
\begin{equation}\label{T0-conclusion1}
\| u_0 \|_{H^1(B_R)} \le C_R \| f\|_{L^2(\mR^d)} \quad \forall \, R > 0. 
\end{equation}
\end{theorem}


We next give an algebraic characterization of  the complementing boundary condition. 

\begin{proposition}\label{pro-complementing}  Let $e$ be a unit vector in $\mR^d$ and let $A_1$ and $A_2$ be two constant positive symmetric matrices. Then  $A_1$ and $A_2$ satisfy the (Cauchy) complementing boundary condition with respect to $e$ if and only if \begin{equation}\label{cond-complementary}
\langle A_2 e, e \rangle \langle A_2 \xi, \xi \rangle  - \langle A_2 e,  \xi \rangle^2 \neq  \langle A_1 e, e \rangle \langle A_1 \xi, \xi \rangle  - \langle A_1 e, \xi \rangle^2 \quad \forall \, \xi \in P \setminus \{0 \}, 
\end{equation}
where
\begin{equation*}
{\cal P} := \big\{\xi \in \mR^d; \langle \xi, e \rangle = 0  \big\}. 
\end{equation*}
In particular, if $A_2 > A_1$ then $A_1$ and $A_2$ satisfy the (Cauchy) complementing boundary condition with respect to $e$. 
\end{proposition}

\begin{remark} \fontfamily{m} \selectfont Assume that $A_1$ is isotropic, i.e., $A_1 = \lambda I$ for some $\lambda > 0$, and $d =2$. Then $A_1$ and $A_2$ satisfy the complementing boundary  condition with respect to  $e$ if and only if $\det A_2 \neq \lambda^2$. In general, \eqref{cond-complementary} is only required  on ${\cal P}$ which is of  co-dimension 1.   
\end{remark}

Using~Theorem~\ref{thm0} and Proposition~\ref{pro-complementing}, one obtains new conditions for which  the well-posedness and  the limiting absorption principle hold. In particular, one can immediately derive  the following  result:

\begin{corollary}\label{cor0} Let $f \in L^2(\mR^d)$ with $\supp f \subset \subset B_{R_0}$, and  
let $u_\delta \in H^1(\mR^d)$ $(0 < \delta < 1)$ be the unique solution of \eqref{Main-eq-delta}. Assume that $A_+: = A  \big|_{\mR^d \setminus \bar D} \in C^1(\bar D_{-\tau})$ and $A_-: = A \big|_{D} \in C^1(\bar D_\tau )$ for some $\tau >0$, and $A_+(x) >  A_-(x)$ or $A_-(x) > A_+(x)$ for all  $x \in \Gamma$.  Then 
\begin{equation*}
\| u_\delta\|_{H^1(B_R)} \le C_R \| f\|_{L^2(\mR^d)} \quad \forall \, R > 0, 
\end{equation*}
for some positive constant $C_R$ independent of $\delta$ and $f$.  Moreover,  $u_\delta \to u_0$ weakly in $H^1_{\loc}(\mR^d)$, as $\delta \to 0$, where  $u_0 \in H^1_{\loc}(\mR^d)$ is the  unique outgoing solution of \eqref{Main-eq}.  We also have 
\begin{equation*}
\| u_0 \|_{H^1(B_R)} \le C_R \| f\|_{L^2(\mR^d)} \quad \forall \, R > 0. 
\end{equation*}
\end{corollary}



To our knowledge, Corollary~\ref{cor0} is new and cannot be obtained using the known approaches mentioned in the introduction.  Corollary~\ref{cor0} is in the same spirit of the one of Bonnet-Ben Dhia, Chesnel, and Ciarlet in \cite{AnneSophieChesnelCiarlet1}; nevertheless, $A_+$ and $A_-$ are not assumed to be isotropic here.

\medskip

The rest of this section contains three subsections. In the first one, we present some  lemmas which are used  in the proof of Theorem~\ref{thm0}. The proof of Theorem~\ref{thm0} is given in the second subsection. 
In the third subsection, we present the proof of Proposition~\ref{pro-complementing}.

\subsection{Preliminaries}

In this section, we present some lemmas used in the proof of Theorem~\ref{thm0}. The first one is on an estimate for solutions to the Helmholtz equation. The proof is based on the unique continuation principle via a compactness argument. 

\begin{lemma}\label{T0-lem1} Let $d \ge 2$, $\Omega$ be a smooth bounded open subset of $\mR^d$,  $f \in L^2(\Omega)$, and let $a$ be a real  uniformly elliptic matrix-valued function and $\sigma$ be a bounded complex function defined in $\Omega$. Assume that $a$ is piecewise Lipschitz and  $v \in H^1(\Omega)$ is a solution to  
\begin{equation*}
\dive(a \nabla  v) + \sigma v = f \mbox{ in } \Omega. 
\end{equation*}
We have 
\begin{equation}
\| v\|_{H^1(\Omega)} \le C \Big( \| f\|_{L^2(\Omega)} + \| v\|_{H^{1/2}(\partial \Omega)} +  \| a \nabla v \cdot \nu \|_{H^{-1/2}(\partial \Omega)} \Big), 
\end{equation}
for some  positive constant $C$ independent of $f$ and $v$. 
\end{lemma}

Here and in what follows, on the boundary  of  a smooth bounded open subset of $\mR^d$,  $\nu$ denotes the normal unit vector directed to its exterior unless otherwise specified. 

\medskip

\noindent{\bf Proof.} We first establish
\begin{equation}\label{T0-lem1-L2}
\| v\|_{L^2(\Omega)} \le C \Big( \| f\|_{L^2(\Omega)} + \| v\|_{H^{1/2}(\partial \Omega)} +  \| a \nabla v \cdot \nu \|_{H^{-1/2}(\partial \Omega)} \Big), 
\end{equation}
by contradiction. Here and in what follows in this proof, $C$ denotes a positive constant independent of $f$, $v$, and $n$. Assume that there exist a sequence $(f_n)  \subset L^2(\Omega)$ and a sequence $(v_n) \subset H^1(\Omega)$ such that 
\begin{equation}\label{T0-lem1-contra1}
\| v_n\|_{L^2(\Omega)} =1, \quad \| f_n\|_{L^2(\Omega)} + \| v_n\|_{H^{1/2}(\partial \Omega)} +  \| a \nabla v_n \cdot \nu \|_{H^{-1/2}(\partial \Omega)} \le 1/n 
\end{equation}
and 
\begin{equation}\label{T0-lem1-contra2}
\dive(a \nabla v_n) +  \sigma v_n = f_n \mbox{ in } \Omega. 
\end{equation}
Multiplying the equation of $\bar v_n$ (the conjugate of $v_n$) and integrating on $\Omega$, we obtain 
\begin{equation}\label{T0-lem1-key}
\| \nabla v_n \|_{L^2(\Omega)} \le C \Big( \| v_n\|_{L^2(\Omega)} + \| f_n\|_{L^2(\Omega)} + \| v_n\|_{H^{1/2}(\partial \Omega)} +  \| a \nabla v_n \cdot \nu \|_{H^{-1/2}(\partial \Omega)} \Big); 
\end{equation}
which implies 
\begin{equation*}
\|v_n \|_{H^1(\Omega)} \le C.  
\end{equation*}
Without loss of generality, one might assume that $v_n \to v$ weakly in $H^1(\Omega)$ and strongly in $L^2(\Omega)$. It follows from \eqref{T0-lem1-contra1} and \eqref{T0-lem1-contra2} that 
\begin{equation*}
\dive(a \nabla v) + \sigma v = 0 \mbox{ in } \Omega
\end{equation*}
and $v = A \nabla v \cdot \nu  = 0 $ on $\partial \Omega$. By the unique continuation principle, see e.g., \cite{Protter60}, $v = 0$ in $\Omega$. This contradicts the fact, by \eqref{T0-lem1-contra1}, 
\begin{equation*}
\| v\|_{L^2(\Omega)} = 1. 
\end{equation*}
Hence \eqref{T0-lem1-L2} holds. The conclusion now follows from \eqref{T0-lem1-key} where $v_n$ is replaced by $v$. \proofend

\begin{remark} \fontfamily{m} \selectfont Assume that $a \in C^1(\bar \Omega)$. Using a three spheres inequality, see e.g., \cite{AlessandriniRondi, MinhLoc2}, one can choose the constant $C$ depending only on $\Omega$, the elliptic and Lipschitz constants of $a$, the boundeness of $a$ and $\sigma$. 
\end{remark}

 The following lemma is used to obtain an a priori estimate for $u_\delta$ defined in \eqref{Main-eq-delta}. 

\begin{lemma} \label{T0-lem2} Let  $f \in L^2(\mR^d)$ with $\supp f \subset \subset B_{R_0}$ and let $u_\delta \in H^1(\mR^d)$ be the unique solution of \eqref{Main-eq-delta}. 
Then 
\begin{equation}\label{est-T0-lem2}
\| u_\delta\|_{H^1(\mR^d)}^2 \le C \Big( \frac{1}{\delta} \Big| \int_{\mR^d} f \bar u_\delta \Big| +  \| f\|_{L^2(B_{R_0})}^2 \Big), 
\end{equation}
for some positive constant $C$ independent of $f$ and $\delta$. Consequently, 
\begin{equation*}
\| u_\delta\|_{H^1(\mR^d)} \le  \frac{C}{\delta}  \| f\|_{L^2(B_{R_0})}. 
\end{equation*}
\end{lemma}

\noindent{\bf Proof.}  Multiplying the equation of $u_\delta$ by $\bar u_\delta$ and integrating on $\mR^d$, we have
\begin{equation}\label{T0-lem2-1}
- \int_{\mR^d} \langle s_\delta A \nabla u_\delta, \nabla u_\delta \rangle +  \int_{\mR^d} k^2 s_0 \Sigma |u_\delta|^2 + i \delta |u_\delta|^2 = \int_{\mR^d} f \bar u_\delta.  
\end{equation}
Considering the imaginary part of \eqref{T0-lem2-1}, we derive that  
\begin{equation*}
\int_{D} |\nabla u_\delta|^2 + \int_{\mR^d} |u_\delta|^2 \le \frac{C}{\delta} \Big|\int_{\mR^d} f \bar u_\delta \Big|. 
\end{equation*}
This implies 
\begin{equation*}
\| u_\delta \|_{H^{1/2}(\partial D)}^2 + \| A \nabla u_\delta \cdot \nu \|_{H^{-1/2}(\partial D)}^2 \le \frac{C}{\delta} \Big|\int_{\mR^d} f \bar u_\delta \Big| +  C \| f\|_{L^2(B_{R_0})}^2. 
\end{equation*}
Let $\Omega$ be the complement of the unbounded connected component of $\mR^d \setminus D$ in $\mR^d$. Applying Lemma~\ref{T0-lem1}, we have 
\begin{equation}\label{T0-lem2-part1}
\| u_\delta\|_{H^1(\Omega)}^2 \le \frac{C}{\delta} \Big|\int_{\mR^d} f \bar u_\delta \Big| + C \| f\|_{L^2(\Omega)}^2. 
\end{equation}
Considering the real part of \eqref{T0-lem2-1} and using \eqref{T0-lem2-part1}, we obtain 
\begin{equation*}
\| u_\delta\|_{H^1(\mR^d)}^2 \le C \Big( \frac{1}{\delta} \Big|\int_{\mR^d} f \bar u_\delta \Big| +  \| f\|_{L^2(B_{R_0})}^2 \Big). 
\end{equation*}
The proof is complete. \proofend


\medskip

The following lemma on the stability of the outgoing solution is standard (see, e.g, \cite{Leis} \footnote{In \cite{Leis}, the proof is given only for $d=2, \, 3$. However, the proof in the case $d > 3$ is similar to the case $d=3$. 
}).

\begin{lemma}\label{lem-Helmholtz} Let  $\Omega \subset B_{R_0}$ be a smooth open subset of $\mR^d$, and let $f \in L^2(\mR^d \setminus \Omega) $ and  $g \in H^{\frac{1}{2}}(\partial \Omega)$. Assume that $\mR^d \setminus \Omega$ is connected, $\supp f \subset B_{R_0}$,   and $v  \in H^1_{\loc}(\mR^d)$ is the unique outgoing  solution of 
\begin{equation*}
\left\{\begin{array}{ll}
\Delta v + k^2  v = f & \mbox{in } \mR^d \setminus \Omega, \\[6pt]
v = g & \mbox{on } \partial \Omega.
\end{array}\right.
\end{equation*}
Then
\begin{equation*}
\| v \|_{H^1(B_r \setminus \Omega)} \le C_r  \big( \| f\|_{L^2(\mR^d \setminus \Omega)} + \| g\|_{H^\frac{1}{2}(\partial \Omega)} \big) \quad \forall \, r > 0,
\end{equation*}
for some positive constants $C_r =C(r, k, \Omega, R_0, d)$.
\end{lemma}

We next recall the following result  \cite[Lemma 2]{Ng-Complementary},  a change of variables formula, which is used repeatedly in this paper.  

\begin{lemma}\label{lem-TO} Let $\Omega_1 \subset \subset \Omega_2 \subset \subset \Omega_3$  be three smooth bounded open subsets of $\mR^d$. Let $a \in [L^\infty(\Omega_2 \setminus \Omega_1)]^{d \times d}$,  $\sigma \in L^\infty(\Omega_2 \setminus \Omega_1)$ and let  $T$ be a diffeomorphism from $\Omega_2 \setminus \bar \Omega_1$ onto $\Omega_3 \setminus \bar \Omega_2$ such that $T(x) = x$ on $\partial \Omega_2$. Assume that   $u \in H^1(\Omega_2 \setminus \Omega_1)$ and set $v = u \circ T^{-1}$. Then
\begin{equation*}
\dive (a \nabla u) +  \sigma u = f \mbox{ in } \Omega_2 \setminus \Omega_1, 
\end{equation*}
for some $f \in L^2(\Omega_2 \setminus \Omega_1)$,  if and only if
\begin{equation}\label{TO-eq}
\dive (T_*a \nabla v) + T_*\sigma v = T_* f \mbox{ in } \Omega_3 \setminus \Omega_2. 
\end{equation}
Moreover, 
\begin{equation}\label{TO-bdry}
v = u \quad \mbox{ and } \quad T_*a \nabla v \cdot \nu = - a \nabla u \cdot \nu  \mbox{ on } \partial \Omega_2.
\end{equation}
\end{lemma}

Recall that $T_*a$, $T_*\sigma$, and $T_*f$ are given in \eqref{TO}.  Here and in what follows, when we mention a diffeomorphism $F:  \Omega \to \Omega'$ for two open subsets $\Omega, \Omega'$ of $\mR^d$, we mean that $F$ is a diffeomorphism, $F \in C^1 (\bar \Omega)$, and $F^{-1} \in C^1(\bar \Omega')$.

\subsection{Proof of Theorem~\ref{thm0}} We first establish  the uniqueness for  \eqref{Main-eq}. Assume that $f$ = 0. We prove that $u_0=0$ if $u_0 \in H^1_{\loc}(\mR^d)$ is an outgoing solution of \eqref{Main-eq}. The proof is quite standard as in the usual case,  in which the coefficients are positive. Multiplying the equation by $\bar u_0$,  integrating on $B_{R}$, and considering the imaginary part, we have, by letting $R \to + \infty$,  
\begin{equation*}
\lim_{R \to + \infty} \int_{\partial B_R} |u_0|^2 = 0. 
\end{equation*}
Here we use the outgoing condition. 
By Rellich's lemma (see,  e.g.,  \cite{ColtonKressInverse}), $u_0 = 0$ in $\mR^d \setminus B_{R_0}$. It follows from the unique continuation principle  that $u_0 = 0$.  The uniqueness is proved. 

\medskip 
We next established \eqref{T0-conclusion0}. Applying Lemma~\ref{T0-lem2}, we have
\begin{equation}\label{T0-part1}
\| u_\delta\|_{H^1(\mR^d)} \le \frac{C}{\delta} \| f \|_{L^2}. 
\end{equation}
In this proof, $C$ denotes a positive constant independent of $\delta$ and $f$. 
Using the difference method due to Nirenberg (see, e.g., \cite{BrAnalyse1}), one has \footnote{We do not claim that $u \in H^2_{\loc}(\mR^d)$; this fact is not true in general.}
\begin{equation}\label{T0-part2}
u_\delta \in H^2(D_{-\tau} \cup D_\tau). 
\end{equation}
For $\tau > 0$ small, define  $F: D_{-\tau} \to D_\tau$ as  follows 
\begin{equation}\label{def-reflection}
F(x_\Gamma + t \nu(x_\Gamma)) = x_\Gamma - t \nu(x_\Gamma) \quad \forall \, x_\Gamma \in \Gamma, \, t \in (-\tau, 0). 
\end{equation}
Let $v_\delta$ be the reflection of $u_\delta$ through $\Gamma$ by  $F$, i.e., 
\begin{equation*}
v_\delta = u_\delta \circ F^{-1} \mbox{ in } D_\tau. 
\end{equation*}
By Lemma~\ref{lem-TO},  we have
\begin{equation*}
\dive(F_*A \nabla v_\delta) + k^2  F_*\Sigma v_\delta +  i \delta  F_* 1 v_\delta = F_*f  \mbox{ in } D_\tau, 
\end{equation*}
and 
\begin{equation*}
v_\delta - u_\delta \big|_{D} = 0, \quad F_*A \nabla v_\delta \cdot \nu - A \nabla u_\delta  \big|_{D} \cdot \nu = i \delta A \nabla u_\delta \cdot \nu \big|_{D} \mbox{ on } \Gamma. 
\end{equation*}
Recall that 
\begin{equation*}
 \dive(A \nabla u_\delta) + k^2  \Sigma u_\delta + k^2 (s_\delta^{-1} s_0 -1)  \Sigma u_\delta + i \delta s_\delta^{-1} u_\delta  = s_\delta^{-1}f  \mbox{ in } D_\tau. 
\end{equation*}
Note that $A_+$ and $A_-$ satisfy the complementing boundary condition on $\Gamma$ if and only if $F_*A_+$ and $A_-$ satisfy the complementing boundary condition on $\Gamma$ by \eqref{cond-complementary} in Proposition~\ref{pro-complementing}.  Applying the result of Agmon, Douglis, and Nirenberg \cite[Theorem 10.2]{ADNII}, we have 
\begin{multline}\label{ADN}
\| u_\delta \|_{H^2(D_{\tau/2})} + \| v_\delta \|_{H^2(D_{\tau/2})}
\\[6pt]
\le C \Big( \| u_\delta \|_{H^1(D_{\tau})} + \| v_\delta \|_{H^1(D_{\tau})} + \| i \delta A \nabla u_\delta  \big|_{D} \cdot \nu \|_{H^{1/2}(\Gamma)} + \|f \|_{L^2} \Big). 
\end{multline}
Since
\begin{equation*}
\|A \nabla u_\delta  \big|_{D}  \cdot \nu  \|_{H^{1/2}(\Gamma)}  \le C \Big( \| u_\delta \|_{H^2(D_{\tau/2})} + \|f \|_{L^2} \Big), 
\end{equation*}
it follows that, for small $\delta$,  
\begin{equation}\label{T0-part3-1}
\| u_\delta \|_{H^2(D_{\tau/2})} + \| v_\delta \|_{H^2(D_{\tau/2})}
\le C \Big(  \| u_\delta \|_{H^1(D_{\tau})} + \| v_\delta \|_{H^1(D_{\tau})} + \| f\|_{L^2} \Big). 
\end{equation}
Using the inequality
\begin{equation*}
\| \varphi \|_{H^1(D_{\tau/2})} \le \lambda \| \varphi\|_{H^2(D_{\tau/2})} + C_\lambda \|\varphi \|_{L^2(D_{\tau/2})}, 
\end{equation*}
we derive from Lemmas~\ref{T0-lem1} and \ref{lem-Helmholtz} that, for small $\delta$,  
\begin{equation}\label{T0-part3}
\| u_\delta \|_{H^2(D_{\tau/2})} + \| v_\delta \|_{H^2(D_{\tau/2})} + \| u_\delta\|_{H^1(B_R)} 
\le C_R \Big( \| u_\delta \|_{L^2(D_{-\tau} \cup D_{\tau})}  + \| f\|_{L^2} \Big) \quad \forall \, R > 0. 
\end{equation}
The proof now follows by a standard compactness argument. We first claim that  
\begin{equation}\label{T0-claim1}
\| u_\delta\|_{L^2(B_{R_0})} \le C \| f\|_{L^2(\mR^d)}. 
\end{equation}
Indeed, assume  that this is not true. By \eqref{T0-part1},  there exist a sequence $(\delta_n) \to 0_+$ and a sequence $(f_n)$ such that $\supp f_n \subset B_{R_0}$, 
\begin{equation*}
\| u_{\delta_n} \|_{L^2(B_{R_0})} = 1, \quad \mbox{ and } \quad \| f_n \|_{L^2(B_{R_0})} \to 0. 
\end{equation*}
We derive from \eqref{T0-part1} and \eqref{T0-part3} that $(u_{\delta_n})$ is bounded in $H^1_{\loc}(\mR^d)$. Without loss of generality, one might assume that $u_{\delta_n} \to u_0$ weakly in $H^1_{\loc}(\mR^d)$ and strongly in $L^2_{\loc}(\mR^d)$ as $n \to + \infty$. Then  $u_0 \in H^1_{\loc}(\mR^d)$, 
\begin{equation*}
\dive(s_0 A \nabla u_0) + k^2 s_0 \Sigma u_0= 0 \mbox{ in } \mR^d, 
\end{equation*}
and $u_0$ satisfies the outgoing condition by the limiting absorption principle. 
It follows that $u_0 = 0$ in $\mR^d$ by the uniqueness. This contradicts the fact $\| u_0\|_{L^2(B_{R_0})} = \lim_{n \to + \infty} \| u_{\delta_n}\|_{L^2(B_{R_0})} = 1$. 
Hence \eqref{T0-claim1} holds. 

A combination of \eqref{T0-part1} and \eqref{T0-claim1} yields 
\begin{equation}\label{T0-part4}
 \| u_\delta\|_{H^1(B_R)} 
\le C_R\| f\|_{L^2}. 
\end{equation}
Hence for any sequence $(\delta_n) \to 0$, there exists a subsequence $(\delta_{n_k})$ such that
$u_{\delta_{n_k}} \to u_0$ weakly in $H^1_{\loc}(\mR^d)$ and strongly in $L^2_{\loc}(\mR^d)$ as $n \to + \infty$. Moreover, $u_0 \in H^1_{\loc}(\mR^d)$, 
\begin{equation*}
\dive(s_0 A \nabla u_0) + k^2 s_0 \Sigma u_0= f \mbox{ in } \mR^d, 
\end{equation*}
and $u_0$ satisfies the outgoing condition. 
Since the limit $u_0$ is unique, $u_\delta \to u_0$ weakly in $H^1_{\loc}(\mR^d)$ and strongly in $L^2(\mR^d)$ as $\delta \to 0$. 
The proof is complete.  \proofend

\subsection{Proof of Proposition~\ref{pro-complementing}}

Using a rotation if necessary, without lost of generality, one may assume that $e = e_d: = (0, \cdots, 0, 1)$. Denote $x = (x', t) \in \mR^{d-1} \times \mR$. Fix a non-zero vector $\xi' = (\xi_1, \cdots, \xi_{d-1}) \in \mR^{d-1}$ and denote $\xi = (\xi', 0)$.  Since $u_j(x) = e^{i \langle x, \xi \rangle} v_j(t)$ ($j=1, 2$)  is a  solution to the equation 
\begin{equation*}
\dive(A_j \nabla u_j ) = 0 \mbox{ in } \mR^{d - 1} \times (0, + \infty),  
\end{equation*}
it follows that, for $j=1, 2$,  
\begin{equation*}
a_j v_j''(t) + 2 i b_j v_j'(t) - c_j v_j(t) = 0 \mbox{ for } t > 0,  
\end{equation*}
where 
\begin{equation*}
a_j  = (A_j)_{d,d}, \quad b_j =  \sum_{k=1}^{d-1} (A_j)_{d, k} \xi_k, \quad \mbox{ and } \quad c_j = \sum_{k=1}^{d-1} \sum_{l =1}^{d-1} (A_j)_{k, l} \xi_k \xi_l.  
\end{equation*}
Here $(A_j)_{k, l}$ denotes the $(k, l)$ component of $A_j$ for $j=1, 2$ and the symmetry of $A_j$ is used. 
Define, for $j=1, 2$,  
\begin{equation*}
\Delta_j = - b_j^2 +  a_j c_j. 
\end{equation*}
Since $A_j$ is symmetric and positive, it is clear that, for $j=1, 2$,  
\begin{equation*}
a_j = \langle A_j e_d, e_d \rangle > 0, \quad b_j = \langle A_j \xi, e_d \rangle,    \quad \mbox{ and }  \quad \Delta_j = \langle A_j e_d, e_d \rangle \langle A_j \xi, \xi \rangle - \langle A_j e_d, \xi \rangle^2 > 0.  
\end{equation*}
Since $v_j$ is required to be bounded, we have 
\begin{equation*}
\quad v_j(t) = \alpha_j e^{\eta_j t}, 
\end{equation*}
where
\begin{equation*}
\eta_j  =( - i b_j - \sqrt{\Delta_j})/ a_j. 
\end{equation*}
Using the fact that $u_1 = u_2 $ and $A_1 \nabla u_1 \cdot e_d = A_2 \nabla u_2 \cdot e_d $, we have 
\begin{equation*}
\alpha_1 = \alpha_2  \quad \mbox{ and } \quad \alpha_1 \Big( \langle i A_2 \xi + \eta_2 A_2 e_d, e_d \rangle - 
 \langle i A_1 \xi + \eta_1 A_1 e_d, e_d \rangle\Big) =0. 
\end{equation*}
The complementing boundary condition is now  equivalent to the fact that 
\begin{equation*}
\Delta_2 \neq \Delta_1, 
\end{equation*}
for all non-zero $ \xi  = (\xi', 0) \in \mR^{d}$.  Condition~\eqref{cond-complementary} is proved. 

It remains to prove that if $A_2 > A_1$ then \eqref{cond-complementary} holds. Define $M = A_2 - A_1$, fix $\xi \in {\cal P} \setminus \{0 \}$,  and set 
\begin{equation*}
\Delta =  \langle A_2 e, e \rangle \langle A_2 \xi, \xi \rangle - \langle A_2 e, \xi \rangle^2 -  \Big( \langle A_1 e, e \rangle \langle A_1 \xi, \xi \rangle - \langle A_1 e, \xi \rangle^2 \Big). 
\end{equation*}
Using the fact $A_2 = A_1 + M$, after a straightforward computation, we obtain
\begin{equation}\label{pro-part0}
\Delta = \langle M e, e \rangle   \langle A_1 \xi, \xi \rangle +  \langle M \xi, \xi \rangle \langle A_1 e, e \rangle +   \langle M e, e \rangle \langle M \xi, \xi \rangle - 2  \langle M e,  \xi  \rangle \langle A_1 e, \xi \rangle - \langle M e,  \xi \rangle^2. 
\end{equation}
We have, by Cauchy's inequality,  
\begin{equation}\label{pro-part1}
\langle M e, e \rangle   \langle A_1 \xi, \xi \rangle  + \langle M \xi, \xi \rangle \langle A_1 e, e \rangle \ge 2 \Big( \langle M \xi, \xi \rangle  \langle M e, e \rangle   \langle A_1 e, e \rangle   \langle A_1 \xi, \xi \rangle \Big)^{1/2}. 
\end{equation}
Since $M$ and $A_1$ are symmetric and positive and $\langle \xi, e \rangle = 0$, we obtain, by Cauchy-Schwarz's inequality, 
\begin{equation}\label{pro-part2}
 \langle M e, e \rangle  \langle M \xi, \xi \rangle   \langle A_1 e, e \rangle   \langle A_1 \xi, \xi \rangle  > \langle M e, \xi \rangle^2 \langle A_1 e,  \xi \rangle^2 . 
\end{equation}
and 
\begin{equation}\label{pro-part3}
 \langle M e, e \rangle \langle M \xi, \xi \rangle > \langle M \xi , e \rangle^2. 
\end{equation}
A combination of \eqref{pro-part0}, \eqref{pro-part1}, \eqref{pro-part2},  and \eqref{pro-part3} yields 
\begin{equation*}
\Delta > 0. 
\end{equation*}
The proof is complete. \proofend

\section{A variational approach via the Dirichlet principle}\label{sect-Thm1}

In this section, we develop a variational method, which is complement to the one in Section~\ref{sect-Thm0}, to deal with a class of $A$ in  which $F_*A_+$ might be equal to $A_-$ on $\Gamma$ and $A_+$ and $A_-$ are not supposed to be smooth near $\Gamma$; this  is not covered by Theorem~\ref{thm0}. One motivation comes from the work of Ola in \cite{Ola}. The other is from the work of  Bonnet-Ben Dhia, Chesnel, and Ciarlet in \cite{AnneSophieChesnelCiarlet1} where the smoothness of the coefficients  is not required. 

\medskip
The following result is the main result of this section. 

\begin{theorem}\label{thm1} Let $f \in L^2(\mR^d)$ with $\supp f \subset B_{R_0}$, and let $u_\delta \in H^1(\mR^d)$ $(0 < \delta < 1)$ be the unique solution of \eqref{Main-eq-delta}.  Assume that there exists a reflection $F$ from $U \setminus D$ onto $D_{\tau}$ for some $\tau>0$ and for  some smooth open set $U \supset \supset D$, i.e., 
$F$ is diffeomorphism and $F(x)  = x$ on $\Gamma$, 
such that 
\begin{equation}\label{T1-cond}
\mbox{either } \quad A - F_*A \ge c \dist(x, \Gamma)^{\alpha} I \quad \mbox{ or  } \quad F_*A - A \ge c \dist(x, \Gamma)^{\alpha} I,  
\end{equation}
on each connected component of $D_\tau$,  for some $c > 0$,   and for some $0 < \alpha < 2$.
Set $v_\delta = u_\delta \circ F^{-1}$ in $D_\tau$. Then 
\begin{equation}\label{T1-conclusion0}
\| u_\delta\|_{L^2(B_R)} + \|u_\delta   - v_\delta\|_{H^1(D_\tau)} + \Big( \int_{D_\tau} \big| \langle (A - F_*A) \nabla u_\delta, \nabla u_\delta \rangle \big| \Big)^{1/2} \le C_R \| f\|_{L^2(\mR^d)}. 
\end{equation}
 Moreover, $u_\delta \to u_0$ weakly  in  $H^1_{\loc}(\mR^d \setminus \Gamma)$ and strongly in $L^2_{\loc}(\mR^d)$ as $\delta \to 0$, where $u_0 \in H^1_{\loc}(\mR^d \setminus \Gamma) \cap L^2_{\loc}(\mR^d)$ is the {\bf unique} outgoing solution of \eqref{Main-eq} such that the LHS of \eqref{T1-stability1} is finite where $v_0 := u_0 \circ F^{-1}$ in $D_\tau$.  Consequently, 
\begin{equation}\label{T1-stability1}
\| u_0 \|_{ L^2(B_{R})} + \| u_0 - v_0 \|_{H^1(D_{\tau})} + \Big( \int_{D_\tau} \big| \langle (A - F_*A) \nabla u_0, \nabla u_0 \rangle \big| \Big)^{1/2} \le C_R \| f\|_{L^2}. 
\end{equation}
Here $C_R$ denotes a positive constant independent of $f$ and $\delta$. 
\end{theorem}

\begin{remark}  \fontfamily{m} \selectfont We only make the assumption on  the lower bound of $F_*A - A$ or $A - F_*A$  in \eqref{T1-cond},  not on the upper bound. 
\end{remark}

The solution $u_0$  in Theorem~\ref{thm1} is not in $H^1_{\loc}(\mR^d)$ as usual.  The meaning of the solution is given in  the following definition. 

\begin{definition} \label{def-solution} \fontfamily{m} \selectfont Let $ f \in L^2(\mR^d)$ with compact support and let $F$ be a reflection  from $U \setminus D$ to $D_{\tau}$ for some $\tau>0$ (small) and for  some open set $D \subset \subset U$, i.e., 
$F$ is diffeomorphism and $F(x)  = x$ on $\Gamma$ such that \eqref{T1-cond} holds. A  function $u_0 \in H^1_{\loc}(\mR^d \setminus \Gamma) \cap L^2_{\loc}(\mR^d)$  such that the LHS of \eqref{T1-stability1} is finite is said to be  a  solution of \eqref{Main-eq} if 
\begin{equation}\label{Def-Weak-eq}
\dive(s_0 A \nabla u_0) + k^2 s_0 \Sigma u_0 =  f \mbox{ in } \mR^d \setminus \Gamma, 
\end{equation}
\begin{equation}\label{Def-Weak-bdry}
u_0 - v_0 = 0 \quad \mbox{ and } \quad (F_*A \nabla v_0 - A \nabla u_0  \big|_{D} ) \cdot \nu = 0 \mbox{ on } \Gamma, 
\end{equation}
and
\begin{equation}\label{Def-Weak-part3}
\lim_{t \to 0_+} \int_{\partial  D_t  \setminus \Gamma} \big( F_*A \nabla v_0 \cdot \nu \bar v_0 - A \nabla u_0 \cdot \nu \bar u_0 \big) = 0. 
\end{equation}
\end{definition}

\begin{remark} \fontfamily{m} \selectfont
Since $u_0 - v_0 \in H^1(D_\tau)$ and  $\dive (F_*A \nabla v_0 - A \nabla u_0) \in L^2(D_{\tau})$ (the LHS of \eqref{T1-stability1} is finite),  it follows that $u_0 -v_0  \in H^{1/2}(\Gamma)$ and $(F_*A \nabla v_0 - A \nabla u_0 \big|_{D} ) \cdot \nu \in H^{-1/2}(\Gamma)$. Hence requirement \eqref{Def-Weak-bdry} makes sense. It is clear that the definition of weak solutions in Definition~\ref{def-solution} coincides with the standard definition of weak solutions  when $\alpha =0$ by Lemma~\ref{lem-TO}. Requirements in \eqref{Def-Weak-bdry} can be seen as  generalized transmission conditions. 
\end{remark}

The proof of Theorem~\ref{thm1} is based on the Dirichlet principle. The  key observation is that the Cauchy data provides the energy of a solution to an elliptic equation (Lemma~\ref{lem1}). The proof is also based on  a new compactness criterion in $L^2$ (Lemma~\ref{lem3}).  The requirement $\alpha < 2$ is used in the compactness argument; we do not know if this condition is necessary. As a direct consequence of Theorem~\ref{thm1} with $\alpha = 0$,  we obtain the following result:


\begin{corollary}\label{cor1} Let $f \in L^2(\mR^d)$ with $\supp f \subset B_{R_0}$, and let $u_\delta \in H^1(\mR^d)$ $(0 < \delta < 1)$ be the unique solution of \eqref{Main-eq-delta}.  Assume that $A \circ F^{-1}(x)$ {\bf or}  $A (x)$ is {\bf isotropic} for every $x \in D_\tau$,  and 
 \begin{equation}\label{Cor1-cond}
\mbox{ either } \quad  A \circ F^{-1}(x) -  A (x) \ge c I  \quad \mbox{ or } \quad A(x) - A \circ F^{-1}(x) \ge c I
\end{equation}
in each connected component $D_\tau$ for some small $\tau > 0$ and for some $c>0$, where $F\big( x_\Gamma + t \nu(x_\Gamma) \big) : = x_\Gamma - t \nu(x_\Gamma)$ for $x_\Gamma \in \Gamma$ and $t \in (-\tau, \tau)$.  Then 
\begin{equation*}
\| u_\delta\|_{H^1(B_R)} \le C_{R} \|f \|_{L^2}. 
\end{equation*}
  Moreover, $u_\delta \to u_0$ weakly  in  $H^1_{\loc}(\mR^d)$ as $\delta \to 0$,  where $u_0 \in H^1_{\loc}(\mR^d)$ is the unique outgoing solution of \eqref{Main-eq} and 
\begin{equation*}
\| u_0\|_{H^1(B_R)} \le C_{R} \|f \|_{L^2}. 
\end{equation*}
\end{corollary}

\begin{remark} \fontfamily{m} \selectfont  Applying  Corollary~\ref{cor1}, one rediscovers  and extends the result obtained by Bonnet-Ben Dhia, Chesnel, and Ciarlet in \cite{AnneSophieChesnelCiarlet1} where $A_+$ and $A_-$ are both isotropic. 
\end{remark}

We next present another consequence of Theorem~\ref{thm1} for the case $\alpha = 1$. The following notation is used. 

\begin{definition}  \fontfamily{m} \selectfont
 The boundary $\Gamma$ of $D$ is called  strictly convex if all its connected components are the boundary of strictly convex sets. 
\end{definition} 

We are ready to present
 
\begin{corollary}\label{cor2} Let $d \ge 3$,  $f \in L^2(\mR^d)$ with $\supp f \subset B_{R_0}$, and let $u_\delta \in H^1(\mR^d)$ $(0 < \delta < 1)$ be the unique solution of \eqref{Main-eq-delta}.  Assume that $D$ is of class $C^3$, $A$ is {\bf isotropic}  and  {\bf constant} in the orthogonal direction of $\Gamma$ in a neighborhood of $\Gamma$, i.e., $A(x_\Gamma + t \nu(x_\Gamma))$ is independent of $t \in (-\tau_0, \tau_0)$  for  $x_\Gamma \in \Gamma$ and for some small positive constant $\tau_0$,  and  $\Gamma$ is  {\bf strictly convex}. There exist $c >0$, $\tau>0$,  a smooth  open set $U \supset \supset D$,  a reflection $F: U \setminus D \to D_\tau$ such that $F_*A - A \ge c \dist(x, \Gamma) I$ or $A - F_*A \ge c \dist(x, \Gamma) I$ on each connected component of $D_\tau$. As a consequence,  $u_\delta$ satisfies \eqref{T1-conclusion0} with $\alpha =1$. 
Moreover, $u_\delta \to u_0$ weakly  in  $H^1_{\loc}(\mR^d \setminus \Gamma)$ as $\delta \to 0$, where $u_0 \in H^1_{\loc}(\mR^d \setminus \Gamma) \cap L^2_{\loc}(\mR^d)$ is the unique outgoing solution of \eqref{Main-eq} and $u_0$ satisfies \eqref{T1-stability1}. 
\end{corollary}

\begin{remark}\fontfamily{m} \selectfont  \label{rem-cor2}  In particular, if $A$ is {\bf isotropic} and {\bf constant}  in each connected component of a neighborhood of $\Gamma$, then the conclusion of Corollary~\ref{cor2} holds. 
\end{remark}

\begin{remark} \fontfamily{m} \selectfont  Applying  Corollary~\ref{cor1}, one rediscovers  and extends the result obtained by Ola \cite{Ola} and Kettunen, Lassas, and Ola  in \cite{KettunenLassas} where $A=I$ in $D$ and $\Gamma$ has one or two connected components. 

\end{remark}

\begin{remark}  \fontfamily{m} \selectfont 
Corollary~\ref{cor2} does {\bf not} hold in two dimensions. The strict convexity of $\Gamma$ is necessary in three dimensions.  In four or higher dimensions, the strict convexity of $\Gamma$ can be relaxed (see Remark~\ref{rem-convexity}). 
\end{remark}

%

The rest of this section containing three subsections is organized as follows. In the first subsection, we present some lemmas used  in the proof of Theorem~\ref{thm1}. The second and the  third subsections are devoted to the proof of Theorem~\ref{thm1} and Corollary \ref{cor2} respectively. 

\subsection{Some useful lemmas}

We begin with the following lemma which plays an important role in the proof of Theorem~\ref{thm1}.

\begin{lemma}\label{lem1} Let  $\Omega$ be a   smooth bounded open subset  of $\mR^d$,  and $A_1$ and $A_2$ be two symmetric uniformly elliptic matrices defined in $\Omega$.  Let $f_1, f_2 \in L^2(\Omega)$, $h \in H^{-1/2}(\partial \Omega)$ and let $u_1, u_2  \in H^1(\Omega)$  be such that 
\begin{equation}\label{L1-eq1}
-\dive(A_1 \nabla u_1) = f_1 \quad \mbox{ and }  \quad - \dive(A_2 \nabla u_2) = f_2 \mbox{ in } \Omega, 
\end{equation}
\begin{equation}\label{L1-bdry1}
u_1 = u_2  \quad \mbox{ and } \quad A_1 \nabla u_1 \cdot \nu = A_2 \nabla u_2 \cdot \nu + h  \mbox{ on } \partial \Omega.  
\end{equation}
Assume that 
\begin{equation}\label{L1-cond}
A_1 \ge A_2  \mbox{ in } \Omega. 
\end{equation}
Then 
\begin{multline}\label{L1-cl1}
\int_{\Omega}  \langle (A_1 - A_2) \nabla u_1, \nabla u_1 \rangle  +  \int_{\Omega}  |\nabla (u_1-u_2)|^2 \\[6pt] 
\le C\Big( \| (f_1, f_2, u_1, u_2) \|_{L^2(\Omega)}^2  + \| h\|_{H^{-1/2}(\partial \Omega)} \| (u_1, u_2)\|_{H^{1/2}(\partial \Omega)} \Big). 
\end{multline}
\end{lemma}

\noindent{\bf Proof.} By considering the real part and the imaginary part separately, without loss of generality, one may assume that all functions in Lemma~\ref{lem1} are real.  Set 
\begin{equation*}
{\cal M }= \| (f_1, f_2, u_1, u_2) \|_{L^2(\Omega)}^2  + \| h\|_{H^{-1/2}(\partial \Omega)} \| (u_1, u_2)\|_{H^{1/2}(\partial \Omega)}. 
\end{equation*}
Multiplying the equation of $u_j$ by $u_j$ (for $j=1, 2$)  and  integrating on $\Omega$,  we have 
\begin{equation}\label{L1-tototo}
\int_{\Omega} \langle A_j \nabla u_j, \nabla u_j \rangle = \int_{\Omega} f_j u_j + \int_{\partial \Omega} A_j \nabla u_j \cdot \nu \,  u_j. 
\end{equation}
Using \eqref{L1-eq1} and \eqref{L1-bdry1}, we derive from \eqref{L1-tototo} that  
\begin{equation}\label{L1-part1*}
\int_{\Omega} \langle A_1 \nabla u_1, \nabla u_1 \rangle - \langle A_2 \nabla u_2, \nabla u_2 \rangle  \le C {\cal M}. 
\end{equation}
Here and in what follows, $C$ denotes a positive constant independent of $f_j$,  $h$,  $u_j$ for $j=1, 2$. 
By the Dirichlet principle, we have
\begin{multline}\label{L1-toto}
\frac{1}{2} \int_{\Omega} \langle A_2 \nabla u_2, \nabla u_2 \rangle - \int_{\Omega} f_2 u_2 - \int_{\partial \Omega} A_2 \nabla u_2 \cdot \nu  \, u_2 \\[6pt]
\le  \frac{1}{2} \int_{\Omega} \langle A_2 \nabla u_1, \nabla u_1 \rangle  - \int_{\Omega} f_2 u_1 - \int_{\partial \Omega} A_2 \nabla u_2 \cdot \nu \,  u_1. 
\end{multline}
A combination of  \eqref{L1-eq1}, \eqref{L1-bdry1},  and \eqref{L1-toto} yields 
\begin{equation}\label{L1-energy1}
 \int_{\Omega} \langle A_2 \nabla u_2, \nabla u_2 \rangle -  \langle A_2 \nabla u_1, \nabla u_1 \rangle 
\le  C {\cal M}. 
\end{equation}
Adding  \eqref{L1-part1*}  and \eqref{L1-energy1}, we obtain
\begin{equation}\label{L1-part0}
\int_{\Omega} \langle (A_1 - A_2) \nabla u_1, \nabla u_1 \rangle \le C {\cal M}.  
\end{equation}
Set 
\begin{equation*}
w = u_1- u_2 \mbox{ in } \Omega. 
\end{equation*}
We have, in $\Omega$,  
\begin{align*}
\dive(A_2 \nabla w) = \dive (A_2 \nabla u_1) - \dive(A_2 \nabla u_2) =  & \dive (A_1 \nabla u_1) - \dive(A_2 \nabla u_2) + \dive([A_2 - A_1]\nabla  u_1)\\[6pt]
= & -f_1 + f_2 +  \dive([A_2 - A_1]\nabla  u_1). 
\end{align*}
Multiplying this equation by $w$, integrating on $\Omega$, we obtain, by \eqref{L1-eq1} and \eqref{L1-bdry1},  
\begin{equation}\label{L1-part1}
\int_{\Omega} |\nabla w|^2 \le \int_{\Omega} C |\langle (A_1 - A_2) \nabla u_1, \nabla w \rangle| + C{\cal M}. 
\end{equation}
Since $A_1 > A_2$ and $A_1$ and $A_2$ are symmetric, we have, for any $\lambda > 0$, 
\begin{equation*}
\int_{\Omega}|\langle (A_1 - A_2) \nabla u_1, \nabla w \rangle|  \le \lambda \int_{\Omega}|\langle (A_1 - A_2) \nabla u_1, \nabla u_1 \rangle| + \frac{1}{4\lambda}  \int_{\Omega}|\langle (A_1 - A_2) \nabla w, \nabla w \rangle|. 
\end{equation*}
It follows from \eqref{L1-part0} and \eqref{L1-part1} that 
\begin{equation}\label{L1-part2}
 \int_{\Omega} |\nabla w|^2 \le C {\cal M}. 
\end{equation}
The conclusion now follows from \eqref{L1-part0} and \eqref{L1-part2}. The proof is complete.   \proofend

\medskip

We next recall Hardy's inequalities (see,  e.g.,  \cite{MazyaSobolev}). 

\begin{lemma}\label{lem-Hardy} Let  $\Omega$ be a smooth bounded open subset  of $\mR^d$. Then, for all $u \in H^1(\Omega)$, and for $\alpha > 1$,  
\begin{equation}\label{Hardy1}
\int_{\Omega} \dist(x, \partial \Omega)^{\alpha - 2} |u(x)|^2 \, dx \le C_{\alpha, \Omega} \int_{\Omega}  \Big( \dist(x, \partial \Omega)^\alpha |\nabla u(x)|^2  + |u(x)|^2 \Big) \, dx. 
\end{equation}
Here $C_{\alpha, \Omega}$ is a positive constant independent of $u$. 
\end{lemma}

\begin{remark}\label{rem-Hardy}  \fontfamily{m} \selectfont  Lemma~\ref{lem-Hardy} also holds for Lipschitz domains, see \cite[Theorem 1.5]{nevcas62}. 

\end{remark}

Using Lemma~\ref{lem-Hardy}, we can prove the following compactness result which is used in the compactness argument in the proof of Theorem~\ref{thm1}.

\begin{lemma}\label{lem3}  Let  $0\le  \alpha < 2$, $\Omega$ be a smooth bounded open subset  of $\mR^d$,  and $(u_n) \subset H^1(\Omega)$. Assume that 
\begin{equation}\label{lem3-cond}
 \sup_{n} \int_{\Omega}  \Big( \dist(x, \partial \Omega)^\alpha |\nabla u_n(x)|^2  + |u_n|^2 \Big) \, dx < + \infty
\end{equation}
Then $(u_n)$ is relatively compact in $L^2(\Omega)$. 
\end{lemma}

\noindent

\noindent{\bf Proof.} Without loss of generality, one can assume that $\alpha > 1$. By Lemma~\ref{lem-Hardy}, we have  
\begin{equation}\label{L3-part1}
\int_{\Omega} \dist(x, \partial \Omega)^{\alpha - 2} |u_n(x)|^2 \, dx \le C_{\alpha, \Omega} \int_{\Omega}  \Big( \dist(x, \partial \Omega)^\alpha |\nabla u_n(x)|^2  + |u_n(x)|^2 \Big) \, dx.
\end{equation}
We derive from  \eqref{lem3-cond} and \eqref{L3-part1} that, for $\tau > 0$ small,  
\begin{equation}\label{L3-part2}
\int_{\Omega_\tau} |u_n(x)|^2 \, dx \le \tau^{2 - \alpha} \int_{\Omega}  \Big( \dist(x, \partial \Omega)^\alpha |\nabla u_n(x)|^2  + |u_n|^2 \Big) \, dx \le C_{\alpha, \Omega} \tau^{2 - \alpha}. 
\end{equation}
Fix $\eps > 0$ arbitrary. Let $\tau > 0$ (small) be such that 
\begin{equation}\label{L3-part3}
\|u_n\|_{L^2(\Omega_\tau)} \le \eps/2 \quad \forall \, n \in \mN. 
\end{equation}
Such a  $\tau$ exists by \eqref{L3-part2}. 
From \eqref{lem3-cond} and Rellich-Kondrachov's compactness criterion, see, e.g., \cite{BrAnalyse1},  there exist $u_{n_1}, \cdots,  u_{n_k}$ such that 
\begin{equation}\label{L3-part4}
\Big\{ u_n \in L^2(\Omega \setminus \Omega_\tau);  \; n \in \mN \Big\} \subset \bigcup_{j=1}^k  \Big\{u \in L^2(\Omega \setminus \Omega_\tau) ; \| u- u_{n_j}\|_{L^2(\Omega \setminus \Omega_\tau)} \le \eps/2 \Big\}. 
\end{equation}
A combination of \eqref{L3-part3} and \eqref{L3-part4} yields 
\begin{equation*}
\big\{ u_n \in L^2(\Omega) \; n \in \mN \big\} \subset \bigcup_{j=1}^k \big\{u \in L^2(\Omega) ; \| u- u_{n_j}\|_{L^2(\Omega)} \le \eps \big\}. 
\end{equation*}
Therefore,  $(u_n)$ is relatively  compact in $L^2(\Omega)$. \proofend

\medskip
We end this section with the following lemma which implies the uniqueness statement in Theorem~\ref{thm1}. 

\begin{lemma}\label{lem-uniqueness} Let   $F$ be a reflection from $U \setminus D$ to $D_{\tau}$ for some small $\tau>0$ and for  some smooth open set $U \supset \supset  D$, i.e., 
$F$ is diffeomorphism and $F(x)  = x$ on $\Gamma$.  
Assume that  $u_0 \in H^1_{\loc}(\mR^d \setminus \Gamma) \cap L^2_{\loc}(\mR^d)$ is an outgoing solution to 
\begin{equation}\label{Eq1}
\dive(s_0 A \nabla u_0) + k^2 s_0  \Sigma u_0 = 0 \mbox{ in } \mR^d \setminus \Gamma, 
\end{equation}
such that the LHS of \eqref{T1-stability1} is finite with $v_0 : = u_0 \circ F^{-1}$ in $D_\tau$, \begin{equation}\label{Def-Weak-bdry*}
u_0 - v_0 = 0 \quad \mbox{ and } \quad (F_*A \nabla v_0 - A \nabla u_0 \big|_{D} ) \cdot \nu = 0 \mbox{ on } \Gamma
\end{equation}
and 
\begin{equation}\label{Def-Weak-part3*}
\lim_{t \to 0_+} \Im\Big\{ \int_{\partial  D_t  \setminus \Gamma} \big( F_*A \nabla v_0 \cdot \nu \bar v_0 - A \nabla u_0 \cdot \nu \bar u_0 \big) \Big\}= 0.  
\end{equation}
Then $u_0=0$ in $\mR^d$. 
\end{lemma}

\noindent{\bf Proof.} Fix $R>R_0$.  Multiplying \eqref{Main-eq} by $\bar u_0$ and integrating on $B_{R} \setminus \big(D \cup F^{-1}(D_{t}) \big)$ and $D \setminus D_{t}$ respectively, one  has, for $0 < t < \tau$,   
\begin{equation}\label{Thm1-part1}
- \int_{B_{R} \setminus \big(D \cup F^{-1}(D_{t}) \big)} \langle A \nabla u_0, \nabla u_0 \rangle + k^2 \int_{B_{R}}  \Sigma |u_0|^2  + \int_{\partial B_R} \partial_r u_0 \bar u_0   +  \int_{\big[\partial F^{-1}(D_t) \big]\setminus \Gamma} A \nabla u_0 \cdot \nu \;  \bar u_0 = 0
\end{equation}
and 
\begin{equation}\label{Thm1-part2}
\int_{D \setminus D_{t}  } \langle A \nabla u_0, \nabla u_0 \rangle - k^2 \int_{D \setminus  D_{t} }  \Sigma |u_0|^2  -  \int_{\partial D_t \setminus \Gamma} A \nabla u_0 \cdot \nu \; \bar u_0  = 0.  
\end{equation}
Here $\nu$ denotes the normal unit vector directed to the exterior of the set in which one integrates. 
Set 
\begin{equation*}
v_0 = u_0 \circ F^{-1} \mbox{ in } D_{\tau}. 
\end{equation*}
Then, by \cite[Lemma 2]{Ng-Complementary},  
\begin{equation*}
\int_{\partial F^{-1}(D_t) \setminus \Gamma} A \nabla u_0 \cdot \nu \;  \bar u_0 = -  \int_{\partial  D_t  \setminus \Gamma} F_*A \nabla v_0 \cdot \nu\;  \bar v_0.
\end{equation*}
It follows from \eqref{Def-Weak-part3*} that 
\begin{equation}\label{Thm1-part3}
\lim_{t \to 0} \left[\int_{\partial F^{-1}(D_t) \setminus \Gamma} A \nabla u_0 \cdot \nu \, \bar u_0 + \int_{\partial  D_t  \setminus \Gamma} A \nabla u_0 \cdot \nu \, \bar u_0 \right]= 0. 
\end{equation}
Adding \eqref{Thm1-part1} and \eqref{Thm1-part2}, letting $t \to 0$, and using \eqref{Thm1-part3}, we obtain
\begin{equation*}
\Im \Big\{ \int_{\partial B_R} \partial_r u_0 \bar u_0 \Big\} = 0. 
\end{equation*}
This implies, by Rellich's lemma, 
\begin{equation*}
u_0 = 0 \mbox{ in } \mR^d \setminus B_{R_0}. 
\end{equation*}
Using \eqref{Def-Weak-bdry*} and  the unique continuation principle,  we reach 
\begin{equation*}
u_0 = 0 \mbox{ in } B_{R_0}. 
\end{equation*}
Hence $u_0 = 0$ in $\mR^d$. The proof is complete. \proofend

\subsection{Proof of Theorem~\ref{thm1}}

The uniqueness of $u_0$ follows from  Lemma~\ref{lem-uniqueness}. 
We next estimate $u_\delta$. 
By Lemma~\ref{T0-lem2}, 
\begin{equation}\label{Thm1-part4}
\| u_\delta\|_{H^1(\mR^d)} \le \frac{C}{\delta} \| f\|_{L^2}. 
\end{equation}
We prove by contradiction that 
\begin{equation}\label{Thm1-contradiction}
\| u_\delta\|_{L^2(B_{R_0})} \le C \| f\|_{L^2}. 
\end{equation}
Suppose that this is not true. There exist $\delta_n \to 0_+$, $f_n \in L^2(\mR^d)$ with $\supp f_n \subset \subset B_{R_0}$ such that 
\begin{equation}\label{Thm1-contradiction1}
\| u_{\delta_n} \|_{L^2(B_{R_0})} = 1 \quad \mbox{ and } \quad \| f_n\|_{L^2} \to 0. 
\end{equation}
Here $u_{\delta_n} \in H^1(\mR^d)$ is the unique solution of \eqref{Main-eq-delta} with $\delta = \delta_n$ and $f = f_n$. Using \eqref{est-T0-lem2} in Lemma~\ref{T0-lem2}, we have 
\begin{equation}\label{Thm1-part4-1}
\| u_\delta\|_{H^1(\mR^d)} \le C\delta_n^{-1/2}. 
\end{equation}
We derive from  Lemma~\ref{lem-TO} that 
\begin{equation*}
\dive(F_*A \nabla v_{\delta_n}) + k^2 F_*\Sigma v_{\delta_n} + i \delta_n F_*1 v_{\delta_n} = F_*f  \mbox{ in } D_\tau, 
\end{equation*}
and 
\begin{equation}\label{Thm1-part5}
v_{\delta_n} = u _{\delta_n} \quad \mbox{ and }  \quad F_*A \nabla v_{\delta_n} \cdot \nu  = (1 + i \delta_n)  A \nabla u_{\delta_n} \Big|_{D}  \cdot \nu \mbox{ on } \Gamma. 
\end{equation}
We also have 
\begin{equation*}
\dive(A \nabla u_{\delta_n}) +  k^2 \Sigma u_\delta +  (s_{\delta}^{-1}s_0 -1) k^2 \Sigma u_{\delta_n} + s_\delta^{-1} i \delta_n u_{\delta_n} = s_\delta^{-1} f \mbox{ in } D_\tau. 
\end{equation*}
Applying Lemma~\ref{lem1} with $D = D_{\tau/2}$ and using \eqref{Thm1-part4-1}, we obtain 
\begin{equation*}
\sup_{n} \int_{D_{\tau/2}} \dist(x, \Gamma)^\alpha |\nabla u_{\delta_n}|^2 + |\nabla (u_{\delta_n} - v_{\delta_n})|^2 < + \infty.
\end{equation*}
By Lemma~\ref{lem3}, 
\begin{equation*}
(u_{\delta_n}), (v_{\delta_n}) \mbox{ are relatively compact in } L^2(D_{\tau/2}). 
\end{equation*}
This implies 
\begin{equation*}
\| \big(  u_{\delta_n}, v_{\delta_n} \big) \|_{H^{1/2}(\partial D_{\tau/4})}, \;  \| \big( A \nabla u_{\delta_n} \cdot \nu, F*A \nabla v_{\delta_n} \cdot \nu \big) \|_{H^{-1/2}(\partial D_{\tau/4})} \mbox{ are bounded.}  
\end{equation*}
From Lemmas~\ref{T0-lem1} and \ref{lem-Helmholtz},  one may assume that   
\begin{equation*}
(u_{\delta_n}) \mbox{ converges in } L^2_{\loc}(\mR^d), 
\end{equation*}
and  $u_{\delta_n}$ and $v_{\delta_n} $ converges almost everywhere.
Let $u_0$ be the limit of $\big(u_{\delta_n}\big)$ in $L^2_{\loc}(\mR^d)$ and $v_0$ be the limit of $\big(v_{\delta_n} \big)$ in $L^2(D_\tau)$. Then $u_0 \in H^1_{\loc}(\mR^d \setminus \Gamma) \cap L^2_{\loc}(\mR^d)$ is a solution to 
\begin{equation*}
\dive (s_0 A \nabla u_0 ) + k^2 s_0 \Sigma u_0 = 0 \mbox{ in } \mR^d \setminus \Gamma, 
\end{equation*} 
$u_0$ satisfies the outgoing condition by the limiting absorption principle, and $v_0 = u_0 \circ F^{-1}$ in $D_\tau$. 
From \eqref{Thm1-part4} and \eqref{Thm1-part5},  we obtain
\begin{equation*}
u_0 - v_0 = 0  \quad \mbox{ and } \quad (A \nabla u_0 \big|_{D} - F_*A \nabla v_0) \cdot \nu = 0 \mbox{ on } \Gamma. 
\end{equation*}
Multiplying the equation of $v_\delta$ and $u_\delta$ by $\bar v_\delta$ and $\bar u_\delta$ respectively, integrating on $D_\tau$, and considering the imaginary part, we have 
\begin{equation*}
\Im \Big\{ \int_{\partial  D_t  \setminus \Gamma} \big( F_*A \nabla v_\delta \cdot \nu \bar v_\delta - (1 + i \delta) A \nabla u_\delta \cdot \nu \bar u_\delta \big)  + \int_{D_\tau} i \delta \langle A \nabla u_\delta, \nabla u_\delta \rangle \Big\} = \Im \Big\{ \int_{D_t} F_*f \bar v_\delta -  \int_{D_t} f \bar u_\delta  \Big\}. 
\end{equation*}
Letting $\delta \to 0$, we obtain 
\begin{equation*}
\Im \Big\{ \int_{\partial  D_t  \setminus \Gamma} \big( F_*A \nabla v_0 \cdot \nu \bar v_0 - A \nabla u_0 \cdot \nu \bar u_0 \big) \Big\} = \Im \Big\{ \int_{D_t} F_*f \bar v_0 -  \int_{D_t} f \bar u_0  \Big\}. 
\end{equation*}
It follows that 
\begin{equation}\label{haha-Thm1}
\lim_{\tau \to 0} \Im \Big\{ \int_{\partial  D_t  \setminus \Gamma} \big( F_*A \nabla v_0 \cdot \nu \bar v_0 - A \nabla u_0 \cdot \nu \bar u_0 \big) \Big\} = 0. 
\end{equation}
Hence $u_0 = 0$ by Lemma~\ref{lem-uniqueness}; this contradicts to the fact $\| u_0\|_{L^2(B_{R_0})} = 1$ by \eqref{Thm1-contradiction1}. Estimate \eqref{Thm1-contradiction} is proved. Estimate~\eqref{T1-conclusion0} now follows from Lemma~\ref{lem1}. Hence, for any sequence $(\delta_n) \to 0$, there exists a subsequence $(\delta_{n_k})$ such that $u_{\delta_{n_k}} \to u_0$ weakly in $L^2_{\loc}(\mR^d)$. It is clear that $u_0 \in H^1_{\loc}(\mR^d) \cap L^2_{\loc}(\mR^d)$, $u_0 - v_0 \in H^1(D_{\tau})$, and $u_0$ is the unique  outgoing condition to 
\begin{equation*}
\dive(s_0A \nabla u_0 ) + k^2 s_0 \Sigma u_0 = f \mbox{ in } \mR^d. 
\end{equation*}
Since the limit $u_0$ is unique, the convergence holds as $\delta \to 0$. It is clear that estimate \eqref{T1-stability1} is a direct  consequence of \eqref{T1-conclusion0}.  The proof is complete. \proofend

\subsection{Proof of Corollary~\ref{cor2}} 

The proof of Corollary~\ref{cor2} is based on a reflection which is different from the standard one used in  Corollary~\ref{cor1}. 
Let $F$ be defined as follows:  
\begin{equation*}
x_\Gamma - t \nu (x_\Gamma) \mapsto x_\Gamma + t \big[1 + t c(x_\Gamma) \big] \nu(x _\Gamma),   
\end{equation*}
for $x_\Gamma \in \Gamma$ and $t > 0$ (small). Here  $c(x_\Gamma) = \beta \tr \Pi(x_\Gamma)$ where $\Pi(x_\Gamma)$ is the second fundamental form of $\Gamma$ at $x_\Gamma$ and $\beta$ is  a constant defined later. In this proof, $\nu(x_\Gamma)$ denotes the unit normal vector of $\Gamma$ at $x_\Gamma$ directed into $D$. Fixing $x_\Gamma \in \Gamma$, we estimate $F_*A - A$ at  $x_\Gamma + t \big[1 + t c(x_\Gamma) \big] \nu(x _\Gamma)$ for small positive $t$. 
To this end, we use local coordinates. Without loss of generality, one may assume that $x_\Gamma = 0$ and around  $x_\Gamma = 0$, $\Gamma$ is presented by the graph of a function $\varphi: (-\eps_0, \eps_0)^{d-1} \to \mR$ with $\varphi(0)  = 0$, and $\{(x', x_d) \in (-\eps_0, \eps_0)^d; x_d > \varphi(x') \} \subset D$. We also assume that $\nabla' \varphi'(0): = (\partial_{x_1}\varphi, \cdots, \partial_{x_{d-1}}\varphi )(0) = 0 \in \mR^{d-1}$ and $\nabla'^2 \varphi(0) = \lambda_1 e_1 \otimes e_1  + \cdots \lambda_{d-1} e_{d-1} \otimes e_{d-1}$ where $\lambda_1, \cdots, \lambda_{d-1}$ are the eigenvalues of $\Pi(x_\Gamma)$.  Here $e_1, \cdots, e_{d}$ is an orthogonal basis of $\mR^d$. Since $\Gamma$ is strictly convex, one can assume that $\varphi$ is strictly convex or strictly concave. We only consider the case  $\varphi$ is  strictly convex;  the other case can be proceeded similarly. Hence,  in what follows, we assume that $\lambda_i > 0$ for $1 \le i \le d-1$.  Set 
\begin{equation*}
\ff(x', t) = \varphi(x', 0).
\end{equation*}
Define
\begin{equation*}
G_1(x', t) = (x', \varphi(x')) + \frac{t\big[1 + t c(x') \big] }{\sqrt{ 1 + |\nabla_{x'}  \varphi(x')|^2}} \Big( -\nabla_{x'} \varphi(x'), 1 \Big). 
\end{equation*} 
A computation yields 
\begin{equation}\label{Cor2-part1}
\nabla G_1(0, t) = I -  t \nabla^2 \ff (0)+ 2 t c(x')  e_d \otimes e_d + O(t^2).   
\end{equation}
Here and in what follows in this paper,  $O(s)$ denotes a quantity or a matrix whose norm is bounded by $C |s|$ for some positive constant $C$ independent of $s$ for small $s$. Define 
\begin{equation*}
G_2(x', t) = (x', \varphi(x')) - \frac{t}{\sqrt{|\nabla' \varphi(x')|^2 + 1}} \Big( -\nabla' \varphi(x'), 1 \Big).
\end{equation*} 
We have
\begin{equation}\label{Cor2-part2}
\nabla G_2(0, t) = I - 2 e_d \otimes e_d + t \nabla^2 \ff(0). 
\end{equation}
From the definition of $F$, $G_1$, and $G_2$,  we have
\begin{equation*}
F(y) = G_1 \circ G_2^{-1}(y). 
\end{equation*}
This yields
\begin{equation*}
\nabla F (y) = \nabla G_1 (x', t) [\nabla G_2 (x', t)]^{-1} \mbox{ where } G_2(x', t) = y. 
\end{equation*}
We derive from \eqref{Cor2-part1} and \eqref{Cor2-part2} that  
\begin{align*}
\nabla F (y) =  I - 2 e_d \otimes e_d - 2t \Pi - 2 t c(0)  \, e_d \otimes e_d + O(t^2),  
\end{align*}
for $y = G_2(0, t)$.  Here for notational ease, we also denote $\Pi = \nabla^2 \ff(0)$. 
We have, for $y = G_2(0, t)$,  
\begin{multline*}
 |\det \nabla F (y)|^{-1} \nabla F (y)^T \nabla F (y)  =   \big[1 + 2t \, \tr \Pi- 2  t c(0) \big] \big(I - 4t \Pi + 4 t c(0)  e_d \otimes e_d \big) + O(t^2) \\[6pt]
=  I + 2 t \sum_{i=1}^{d-1} \big[\tr \Pi  - 2\lambda_i -    c(0) \big] e_i \otimes e_i + 2 t \big[ c (0) + \tr \Pi] e_d \otimes e_d + O(t^2). 
\end{multline*}
By taking $c(0) = \beta \tr \Pi$ with  $-1< \beta < 0$ and $\beta$ is closed to $-1$, we have 
\begin{equation*}
B: =  |\det \nabla F (y)|^{-1} \nabla F (y)^T \nabla F (y) - I   \ge  \gamma   t I. 
\end{equation*}
The conclusion now follows from Theorem~\ref{thm1}. The proof is complete. \proofend

\begin{remark} \label{rem-convexity} \fontfamily{m} \selectfont Corollary~\ref{cor2} does not hold for $d=2$.  Indeed, assume that  $A = I$ in $\mR^2$, $D = B_{r_2} \setminus B_{r_1}$ for $0 < r_1 < r_2$. Let $F: B_{r_2^2/ r_1} \setminus B_{r_2} \to B_{r_2} \setminus B_{r_1}$ be the Kelvin transform with respect to $\partial B_{r_2}$ and let $\Sigma = F_*1$ in $B_{r_2} \setminus B_{r_1}$, then $F_*A = A$ and $F_*\Sigma = \Sigma$: the resonance appears (Proposition~\ref{pro-resonant} in Section~\ref{sect-resonant}).  The strict convexity condition of $\Gamma$ is necessary in three dimensions. In fact assume that  $D = \{(x_1, x_2, x_3) \in \mR^3; \; x_1^2 + x_2 ^2 < 1 \mbox{ and } 0 < x_3 < 1 \}$ and let $G: \mR^2 \times (0, 1) \setminus D \to D$ be defined by $G(x_1, x_2, x_3) = \big(F(x_1, x_2), x_3 \big)$. Set $(A, \Sigma) = (I, 1)$ in $\mR^3 \setminus D$ and $(I, G_*1)$ otherwise. The problem is not well-posed again for some $f$ by Proposition~\ref{pro-resonant} in Section~\ref{sect-resonant}. 
Nevertheless, the strict convexity condition can be weaken in four or higher dimensions. To illustrate this point, let consider the case $d=4$. Then 
\begin{align*}
\frac{1}{2t} B  = &  (\lambda_2 + \lambda_3 - \lambda_1 - \beta) e_1 \otimes e_1 + (\lambda_1 + \lambda_3 - \lambda_2 - \beta) e_e \otimes e_2 \\[6pt]
&   + (\lambda_1 + \lambda_2 - \lambda_2 - \beta) e_3 \otimes e_3 +  (1 + \beta) (\lambda_1 + \lambda_2 + \lambda_3) e_4 \otimes e_4 + O(t).
\end{align*}
Assume that  $\lambda_1, \lambda_2, \lambda_3 \ge 0$ and if $\lambda_1 \lambda_2 \lambda_3 = 0$ then only one of them is 0.  Then 
$B \ge \gamma t I$ if $\beta$ is chosen as in the proof of Corollary~\ref{cor2}. Hence the conclusion of Corollary~\ref{cor2} holds in this case. 
\end{remark}

\section{A variational approach via the multiplier technique}\label{sect-Thm2}

In this section we develop a variational approach via the multiplier technique to deal with the case $F_*A = A$ in $D_\tau$. This  is complement to the results  in the previous sections. The main result of this section  is:

\begin{theorem}\label{thm2} Let $f \in L^2(\mR^d)$ with $\supp f \subset B_{R_0} \setminus \Gamma$, and let $u_\delta \in H^1(\mR^d)$ $(0 < \delta < 1)$ be the unique solution of \eqref{Main-eq-delta}.  Assume that there exists a reflection $F$ from $U \setminus D$ to $D_{\tau}$, for some $\tau>0$ and for some smooth bounded open subset $U \supset \supset D$,  i.e., 
$F$ is diffeomorphism and $F(x)  = x$ on $\partial D$, 
such that  either
\begin{equation}\label{T2-cond1}
F_*A - A \ge 0 \quad \mbox{ and } \quad \Sigma - F_*\Sigma  \ge c \dist(x, \Gamma)^\beta,   
\end{equation}
or
\begin{equation}\label{T2-cond2}
A - F_*A \ge 0 \quad  \mbox{ and } \quad  F_*\Sigma - \Sigma   \ge c \dist(x, \Gamma)^\beta,    
\end{equation}
in each connected component of $D_\tau$, for some $\beta > 0$ and  $c > 0$.  Set $v_\delta = u_\delta \circ F^{-1}$ in $D_\tau$. Then, for all $0< \rho < R$,  
\begin{multline}\label{T2-conclusion0}
\int_{B_{R} \setminus (D_\rho \cup D_{-\rho })} |u_\delta|^2 +  \int_{D_\tau} |\Sigma - F_*\Sigma| |u_\delta|^2  +  \int_{D_\tau}  | \langle (A - F_*A) \nabla u_\delta, \nabla u_\delta \rangle |\\[6pt] 
+ \int_{D_\tau} |u_\delta   - v_\delta |^2 + |\nabla (u_\delta - v_\delta)|^2  \le C_{R, \rho} \| f\|_{L^2(\mR^d)}^2. 
\end{multline}
Moreover, $u_\delta \to u_0$ weakly  in  $H^1_{\loc}(\mR^d \setminus \Gamma)$ and strongly in $L^2_{\loc}(\mR^d \setminus \Gamma)$ as $\delta \to 0$, where $u_0 \in H^1_{\loc}(\mR^d \setminus \Gamma)$ is the {\bf unique} outgoing solution of \eqref{Main-eq} 
such that the LHS of \eqref{T2-stability1} is finite where $v_0 := u_0 \circ F^{-1}$ in $D_\tau$. Consequently, 
\begin{multline}\label{T2-stability1}
\int_{B_{R} \setminus (D_\rho \cup D_{-\rho })} |u_0|^2 +  \int_{D_\tau} |\Sigma - F_*\Sigma| |u_0|^2  +  \int_{D_\tau}  | \langle (A - F_*A) \nabla u_0, \nabla u_0 \rangle |\\[6pt] 
+ \int_{D_\tau} |u_0   - v_0|^2 + |\nabla (u_0- v_0)|^2  \le C_{R, \rho} \| f\|_{L^2(\mR^d)}^2. 
\end{multline}
Here $C_{R, \rho}$ denotes a positive constant depending on $R$, $\rho$,  $A$, $\Sigma$, $R_0$, $\beta$, $c$, and the distance between $\supp f$ and $\Gamma$,   but independent of $f$ and $\delta$. 
\end{theorem}

The solution $u_0$  in Theorem~\ref{thm2} is not in $L^2_{\loc}(\mR^d)$.  Its meaning is given in  the following definition. 

\begin{definition} \label{def-2} \fontfamily{m} \selectfont Let $ f \in L^2(\mR^d)$ with $\supp f \subset \subset \mR^d \setminus \Gamma$,  and let $F$ be a reflection from $U \setminus D$ to $D_{\tau}$ for some $\tau>0$  and for  some smooth open set $U \supset \supset D$ such that  \eqref{T2-cond1} or \eqref{T2-cond2} holds.  A  function $u_0 \in H^1_{\loc}(\mR^d \setminus \Gamma)$ such that the LHS of \eqref{T2-stability1} is finite is called a  solution to \eqref{Main-eq} with $v_0 = u_0 \circ F^{-1}$ if 
\begin{equation}\label{Def-Weak-eq1}
\dive(s_0 A \nabla u_0) + k^2 s_0 \Sigma u_0 =  f \mbox{ in } \mR^d \setminus \Gamma, 
\end{equation}
\begin{equation}\label{Def-Weak-bdry1}
u_0 - v_0 = 0 \quad \mbox{ and } \quad (F_*A \nabla v_0 - A \nabla u_0 \big|_{D} ) \cdot \nu = 0 \mbox{ on } \Gamma, 
\end{equation}
and 
\begin{equation}\label{Def-Weak-part3-1}
\lim_{t \to 0_+} \Im\Big\{ \int_{\partial  D_t  \setminus \Gamma} \big( F_*A \nabla v_0 \cdot \nu \bar v_0 - A \nabla u_0 \cdot \nu \bar u_0 \big) \Big\}= 0.  
\end{equation}

\end{definition}

\begin{remark} \fontfamily{m} \selectfont
Since $u_0 - v_0 \in H^1(D_{\tau})$, $(\Sigma - F_* \Sigma ) u_0 \in L^2(D_\tau)$,  and $(A - F_* A) \nabla u_0 \in L^2(D_\tau)$,  it follows that   $\dive (F_*A \nabla v_0 - A \nabla u_0) \in L^2(D_{\tau})$ and $F_*A \nabla v_0 - A \nabla u_0 \in L^2(D_\tau)$. Hence  $u_0 -v_0  \in H^{1/2}(\Gamma)$,  and $(F_*A \nabla v_0 - A \nabla u_0 \big|_{D}) \cdot \nu \in H^{-1/2}(\Gamma)$. Requirement \eqref{Def-Weak-bdry} makes sense. 
\end{remark}

\begin{remark}  \fontfamily{m} \selectfont $\beta$ is only required to be positive in Theorem~\ref{thm2}. 
In \eqref{T2-cond1} and \eqref{T2-cond2} of Theorem~\ref{thm2}, we only make the assumption on  the lower bound and  not on the upper bound of the quantities considered. 
\end{remark}

The proof of Theorem~\ref{thm2} is based on a variational approach via the multiplier technique. One of the key point of the proof is Lemma~\ref{lem5}, a variant of Lemma~\ref{lem1},  where the multipliers are used.  Sylvester in \cite{Sylvester} used related ideas to study the transmission eigenvalues problem. The compactness argument used in the proof of Theorem~\ref{thm2} is different from the standard one used in the proof of Theorem~\ref{thm1} due to fact the family $(u_\delta)$ might not be bounded in $L^2_{\loc}(\mR^d)$ in this context considered in Theorem~\ref{thm2}. 

\medskip 
Here is a corollary of Theorem~\ref{thm2} which is a complement to   Corollary~\ref{cor2} in two dimensions.  

\begin{corollary}\label{cor3} Let $d =2$, $\sigma_0 \in \mR$, $D = B_1$,   $f \in L^2(\mR^d)$ with $\supp f \subset B_{R_0} \setminus \Gamma$, and let $u_\delta \in H^1(\mR^d)$ $(0 < \delta < 1)$ be the unique solution of \eqref{Main-eq-delta}.  Assume that $(A, \Sigma) = (I, 1)$ in $D_{-\tau}$ and $(A, \Sigma) = (I, \sigma_0)$ in $D_{\tau}$ for some small $\tau > 0$. Let  $F$ be the Kelvin transform with respect to $\partial D$. Then $F_*I = I$ in $D_{\tau/2}$ and $F_*1 -  \sigma_0 \ge c \dist(x, \Gamma) $ in  $D_{\tau/2}$ if $\sigma_0 \le 1$ and $\sigma_0- F_*1 > c $  in  $D_{\tau/2}$ if $\sigma_0 > 1$ for some $c>0$.  As a consequence, $u_\delta \to u_0$ weakly  in  $H^1_{\loc}(\mR^2 \setminus \Gamma)$ as $\delta \to 0$, where $u_0 \in H^1_{\loc}(\mR^d \setminus \Gamma)$ is the unique outgoing solution of \eqref{Main-eq}; moreover, $u_0$ satisfies \eqref{T2-stability1}. 
\end{corollary}

The rest of this section contains three subsections and is devoted to the proof of Theorem~\ref{thm2} and Corollary~\ref{cor3}. The first one is on a variant of Lemma~\ref{lem1} used in the proof of Theorem~\ref{thm2}. The proof of Theorem~\ref{thm2} and Corollary~\ref{cor3} are given in the last two subsections. 

\subsection{A useful lemma}

The following lemma is a variant of Lemma~\ref{lem1} and plays an important role in the proof of Theorem~\ref{thm2}. 

\begin{lemma}\label{lem5} Let $\Omega$ be a   smooth bounded open subset  of $\mR^d$,  and $A_1$ and $A_2$ be two uniformly elliptic matrices, and $\Sigma_1$ and $\Sigma_2$ be two bounded real functions defined in $\Omega$. Let  $f_1, f_2 \in L^2(\Omega)$, $h \in H^{-1/2}(\partial \Omega)$,  and let $u_1, u_2 \in H^1(\Omega)$  be such that
\begin{equation}\label{lem5-eq*}
\dive(A_1 \nabla u_1) +  \Sigma_1 u_1  = f_1 \quad \mbox{ and }  \quad  \dive(A_2 \nabla u_2) +  \Sigma_2 u_2= f_2 \mbox{ in } \Omega, 
\end{equation}
and 
\begin{equation}\label{lem5-bdry*}
u_1 = u_2  \quad \mbox{ and } \quad A_1 \nabla u_1 \cdot \nu = A_2 \nabla u_2 \cdot \nu + h  \mbox{ on } \partial \Omega.  
\end{equation}
Assume that 
\begin{equation}\label{lem5-cond}
A_1 \ge  A_2    \quad \mbox{ and } \quad \Sigma_2 \ge  \Sigma_1   \mbox{ in } \Omega, 
\end{equation}
We have
\begin{equation*}
\int_{\Omega} (\Sigma_2 - \Sigma_1) |u_2|^2 + \langle (A_1 - A_2) \nabla u_2, \nabla u_2 \rangle + |\nabla (u_1 - u_2)|^2  \le C {\cal N} (f_1, f_2, h, u_1, u_2), 
\end{equation*}  
where, for some positive constant $C$ independent of $u_1$, $u_2$, $f_1$, $f_2$, and $h$,  
\begin{multline*}
{\cal N}(f_1, f_2, h, u_1, u_2) = \|(u_1, u_2) \|_{L^2(\Omega)}\|(f_1, f_2) \|_{L^2(\Omega)}\\[6pt] +  \|h\|_{H^{-1/2}(\partial \Omega)} \|(u_1, u_2) \|_{H^{1/2}(\partial \Omega)} + \|u_1 - u_2 \|_{L^2(\Omega)}^2. 
\end{multline*}
\end{lemma}


\noindent{\bf Proof.} By considering the real part and the imaginary part separately, without loss of generality, one may assume that all functions mentioned in Lemma~\ref{lem5} are real.  Define 
\begin{equation*}
w = u_1 - u_2 \mbox{ in } \Omega. 
\end{equation*}
From \eqref{lem5-eq*} and \eqref{lem5-bdry*}, we have
\begin{equation}\label{lem5-eq}
\dive(A_1 \nabla w) +  \Sigma_1 w = f_1 - f_2 + (\Sigma_2 - \Sigma_1) u_2  + \dive ([A_2 - A_1] \nabla u_2) \mbox{ in } \Omega, 
\end{equation}
\begin{equation}\label{lem5-bdry}
w = 0, \quad \mbox{ and } \quad A \nabla w \cdot \nu = h \mbox{ on } \partial \Omega. 
\end{equation}
Multiplying  \eqref{lem5-eq} by $u_2$ and  integrating on $\Omega$,  we have 
\begin{multline*}
\int_{\Omega}( f_1 - f_2) u_2 + (\Sigma_2 - \Sigma_1) |u_2|^2 + \int_{\Omega} \langle (A_1 - A_2) \nabla u_2, \nabla u_2 \rangle\\[6pt]
=  \int_{\Omega} \Big(\dive(A_1 \nabla w) +  \Sigma_1 w \Big) u_2 +  \int_{\partial \Omega} (A_1 - A_2) \nabla u_2 \cdot \nu  \; u_2. 
\end{multline*}
Integrating by parts and using the fact that
\begin{equation*}
A_1 \nabla w \cdot \nu + (A_1 - A_2) \nabla u_2  \cdot \nu= A_1 \nabla u_1 \cdot \nu - A_2 \nabla u_2 \cdot \nu = h \mbox{ on } \partial \Omega, 
\end{equation*}
\begin{equation*}
\int_{\Omega} A_1 \nabla w \nabla u_2 = \int_{\Omega} A_2 \nabla w \nabla u_2 + \int_{\Omega} (A_1-A_2) \nabla w \nabla u_2, 
\end{equation*}
and
\begin{equation*}
\Sigma_1 w u_2 = (\Sigma_1 - \Sigma_2)w u_2 + \Sigma_2 w u_2, 
\end{equation*}
we derive from \eqref{lem5-cond} and  \eqref{lem5-bdry} that 
\begin{multline}\label{lem5-part1}
\int_{\Omega} (\Sigma_2 - \Sigma_1) |u_2|^2 + \int_{\Omega} \langle (A_1 - A_2) \nabla u_2, \nabla u_2 \rangle \\[6pt] 
\le    C{\cal N} (f_1, f_2, h, u_1,u_2) + \int_{\Omega} (\Sigma_1 - \Sigma_2) w u_2 + \int_{\Omega}\langle (A_2 - A_1) \nabla w, \nabla u_2 \rangle. 
\end{multline}
Here and in what follows in this proof, $C$ denotes a positive constant independent of $f_j$, $h$, $u_j$ for $j=1, 2$. 
Multiplying \eqref{lem5-eq} by $w$ and integrating on $\Omega$, we have
\begin{equation}\label{lem5-part2}
\int_{\Omega} \langle A_1 \nabla w, \nabla w \rangle  \le C {\cal N}(f_1, f_2, h, u_1, u_2)  +   \int_{\Omega} (\Sigma_1 - \Sigma_2) w u_2 +  \int_{\Omega} \langle (A_2 - A_1) \nabla u_2, \nabla w \rangle. 
\end{equation}
A combination of \eqref{lem5-part1} and \eqref{lem5-part2} yields 
\begin{multline}\label{lem5-part2*}
\int_{\Omega} (\Sigma_2 - \Sigma_1) |u_2|^2 + \int_{\Omega} \langle (A_1 - A_2) \nabla u_2, \nabla u_2 \rangle + \int_{\Omega} \langle A_1 \nabla w, \nabla w \rangle  \\[6pt] 
\le    C {\cal N} (f_1, f_2, h, u_1,u_2) +  2 \int_{\Omega}  (\Sigma_1 - \Sigma_2) w u_2 +  2 \int_{\Omega}\langle (A_2 - A_1) \nabla w, \nabla u_2 \rangle . 
\end{multline}
We have, since $A_1 \ge A_2$, 
\begin{equation}\label{lem5-t2}
2 \int_{\Omega}  \langle (A_2 - A_1) \nabla u_2, \nabla w  \rangle  \le \lambda \int_{\Omega} \langle (A_1 - A_2) \nabla u_2, \nabla u_2 \rangle  + \lambda^{-1} \langle (A_1 - A_2) \nabla w, \nabla w \rangle, 
\end{equation} 
and, since $\Sigma_2 \ge \Sigma_1$,  
\begin{equation}\label{lem5-t3}
 2 \int_{\Omega}  (\Sigma_1 - \Sigma_2) w u_2  \le 2 \int_{\Omega}  (\Sigma_2 - \Sigma_1) w^2 + \frac{1}{2}  \int_{\Omega}  (\Sigma_2 - \Sigma_1)u_2^2. 
 \end{equation}
By choosing $\lambda$ smaller than 1 and close to 1, we derive from \eqref{lem5-part2*}, \eqref{lem5-t2},  and \eqref{lem5-t3} that 
\begin{equation}\label{lem5-part2**}
\int_{\Omega} (\Sigma_2 - \Sigma_1) |u_2|^2 + \int_{\Omega} \langle (A_1 - A_2) \nabla u_2, \nabla u_2 \rangle + \int_{\Omega} |\nabla w|^2 
\le    C{\cal N} (f_1, f_2, h, u_1,u_2). 
\end{equation}
The proof is complete. \proofend


\subsection{Proof of Theorem~\ref{thm2}} 

The proof of  the uniqueness of $u_0$, i.e., if $f = 0$ then $u_0=0$ is similar to the one of Lemma~\ref{lem-uniqueness}. The details are left to the reader. 

\medskip
We next establish the estimate for $u_\delta$ by a compactness argument. The compactness argument used in this proof  is different from the one in the proof of Theorem~\ref{thm1} due to the loss of the control of $u_\delta$ in $L^2_{\loc}(\mR^d)$.   Without loss of generality, one may assume that $\supp f \cap (D_{\tau} \cup F^{-1}(D_\tau)) = \O$.  By Lemma~\ref{T0-lem2},  we have
\begin{equation}\label{Thm2-important}
\| u_\delta \|_{H^1(\mR^d)}^2 \le \frac{C}{\delta}  \| f\|_{L^2(B_{R_0})} \| u_\delta\|_{L^2(B_{R_0} \setminus (D_\tau \cup F^{-1} (D_\tau)))}. 
\end{equation}
We first prove  that
\begin{equation}\label{Thm2-claim1}
 \| u_\delta\|_{L^2(B_{R_0} \setminus (D_{\tau_1} \cup F^{-1} (D_{\tau_1})))}
 \le  C \| f\|_{L^2}, 
\end{equation}
by contradiction \footnote{We do not prove that $ \| u_\delta\|_{L^2(B_{R_0})} \le  C \| f\|_{L^2}$. This is different from the proof of Theorem~\ref{thm1}.} where $0 < \tau_1 < \tau/3$ is a positive constant chosen later.  Assume that there exist $\delta_n \to 0$, $f_n \in L^2(\mR^d)$ with $\supp f_n \subset B_{R_0}$ and $\supp f_n \cap (D_{\tau} \cup F^{-1}(D_{\tau}) ) = \O$ such that 
\begin{equation}\label{Thm2-contradiction}
 \| f_n\|_{L^2} \to 0 \quad \mbox{ and } \quad \| u_n\|_{L^2(B_{R_0} \setminus (D_{\tau_1} \cup F^{-1} (D_{\tau_1})))} =1, 
\end{equation}  
where $u_n$ is the solution of \eqref{Main-eq-delta} with $\delta = \delta_n$ and $f = f_n$. 
Set $v_n = u_n \circ F^{-1} \mbox{ in } D_\tau$. 
By Lemma~\ref{lem-TO},  
\begin{equation}\label{Thm2-part1}
\dive(F_*A \nabla v_n) + k^2 F_* \Sigma v_n + i \delta_n F_*1 v_n  = 0 \mbox{ in }  D_\tau,  
\end{equation}
and
\begin{equation}\label{Thm2-transmission}
v_n = u_n \quad \mbox{ and } \quad A \nabla v_n \big|_{D} \cdot \nu = F_*A \nabla u_n \cdot \nu + i \delta_n A \nabla u_n \big|_{D}\cdot \nu \mbox{ on } \Gamma. 
\end{equation}
We also have 
\begin{equation}\label{Thm2-part2}
\dive(A \nabla u_n) + k^2 \Sigma u_n +  \Big( i \delta_n s_{\delta_n}^{-1} +  [s_0s_{\delta_n}^{-1}-1] k^2 \Sigma \Big)u_n  = 0 \mbox{ in } D_\tau. 
\end{equation}
Applying Lemma~\ref{lem5} with $D = D_{\tau/2}$ and using \eqref{Thm2-important} and \eqref{Thm2-contradiction}, we have 
\begin{multline}\label{Thm2-part2**}
 \int_{D_{\tau/2}} |\Sigma - F_*\Sigma| |u_n|^2  +  \int_{D_{\tau/2}}  | \langle (A - F_*A) \nabla u_n, \nabla u_n \rangle |\\[6pt] 
+ \int_{D_{\tau/2}} |\nabla (u_n - v_n)|^2   \le C_\tau \Big(1 +  \int_{D_{\tau_1}} |u_n - v_n|^2 \Big). 
\end{multline}
By choosing $\tau_1$ small enough, one has 
\begin{equation*}
C_\tau \int_{D_{\tau_1}} |u_n - v_n|^2  \le \frac{1}{2} \int_{D_{\tau/2}} |\nabla (u_n - v_n)|^2, 
\end{equation*}
since $u_n - v_n = 0$ on $\Gamma$.  It follows from \eqref{Thm2-part2**} that 
\begin{equation}\label{Thm2-part2*}
 \int_{D_{\tau/2}} |\Sigma - F_*\Sigma| |u_n|^2  +  \int_{D_{\tau/2}}  | \langle (A - F_*A) \nabla u_n, \nabla u_n \rangle |
+ \int_{D_{\tau/2}} |\nabla (u_n - v_n)|^2 +  \int_{D_{\tau/2}}  |u_n - v_n|^2   \le C_\tau. 
\end{equation}
This implies, by \eqref{T2-cond1} and \eqref{T2-cond2}, for $0 < \rho < \tau/4$, 
\begin{equation*}
\| \big(  u_n, v_n \big) \|_{H^{1/2}(\partial D_{\rho})}, \;  \| \big( A \nabla u_n \cdot \nu, F*A \nabla v_n \cdot \nu \big) \|_{H^{-1/2}(\partial D_{\rho})} \mbox{ are bounded.}  
\end{equation*}
Using Lemmas~\ref{T0-lem1} and \ref{lem-Helmholtz}, we derive  that 
\begin{equation}\label{Thm2-part2*-1}
\int_{B_R \setminus (D_\rho \cup D_\rho)} |u_n|^2 + |\nabla u_n|^2  \le C_{\rho, R}, 
\end{equation}
for $0 < \rho < R$. 
Without loss of generality, one may assume that $u_n \to u_0$ weakly in $H^1_{\loc}(\mR^d \setminus \Gamma)$, and strongly in $L^2_{\loc}(\mR^d \setminus \Gamma)$, $v_n \to v_0$ weakly in $H^1_{\loc} (D_\tau)$, and $u_n  - v_n \to u_0 - v_0$ weakly in $H^1(D_\tau)$ and strongly in $L^2_{\loc}(D_\tau)$, and $v_0 = u_0 \circ F^{-1}$ in $D_\tau$. We have, by \eqref{Thm2-part2*}, 
 \begin{multline*}
 \int_{D_\tau} |\Sigma - F_*\Sigma| |u_0|^2  +  \int_{D_\tau}  | \langle (A - F_*A) \nabla u_0, \nabla u_0 \rangle |\\[6pt] 
+ \int_{D_\tau} |u_0   - v_0 |^2 + |\nabla (u_0 - v_0)|^2 + \int_{B_R \setminus (D_\rho \cup D_\rho)} |u_0|^2 + |\nabla u_0|^2  \le C_{\rho, R}, 
\end{multline*}
and $u_0 \in H^1_{\loc}(\mR^d \setminus \Gamma)$ is an outgoing solution to the equation 
\begin{equation*}
\dive (s_0 A \nabla u_0) + k^2 s_0 \Sigma u_0 = 0 \mbox{ in } \mR^d \setminus \Gamma. 
\end{equation*}
From \eqref{Thm2-important} and \eqref{Thm2-transmission}, we obtain 
\begin{equation*}
u_0 - v_0 = 0 \quad \mbox{ and } (A \nabla u_0\big|_{D} - F_*A \nabla v_0) \cdot \nu = 0 \mbox{ on } \Gamma. 
\end{equation*}
Similar to \eqref{haha-Thm1}, we also have
\begin{equation*}
\lim_{t \to 0_+} \Im\Big\{ \int_{\partial  D_t  \setminus \Gamma} \big( F_*A \nabla v_0 \cdot \nu \bar v_0 - A \nabla u_0 \cdot \nu \bar u_0 \big) \Big\}= 0.  
\end{equation*}
Hence $u_0 = 0$ in $\mR^d$ by the uniqueness. We have a contradiction with the fact that \\ $\| u_0\|_{L^2(B_{R_0} \setminus (D_{\tau_1} \cup F^{-1} (D_{\tau_1})))} = \lim_{n \to \infty} \| u_n\|_{L^2(B_{R_0} \setminus (D_{\tau_1} \cup F^{-1} (D_{\tau_1})))} = 1$. Claim \eqref{Thm2-claim1} is proved. 
The conclusion now is standard as in the proof of Theorem~\ref{thm1}. The details are left to the reader.  \proofend

\subsection{Proof of Corollary~\ref{cor3}}

It suffices  to check $F_* 1 - \sigma_0 \ge c \dist(x, \Gamma)$ if $ \sigma_0 \le 1$ and $\sigma_0 - F_* 1 > c$ if $\sigma_0 > 1$ in $D_{\tau/2}$ for some $c > 0$ provided that  $\tau$ is small enough.  A computation gives
\begin{equation*}
\det (\nabla F)(y) = 1 - 4 \dist(x, \Gamma) + O\big( \dist(x, \Gamma)^2\big), 
\end{equation*}
where $F(y) = x$. 
This implies 
\begin{equation*}
1/ \det(\nabla F) (y) = 1 + 4 \dist(x, \Gamma) + O\big( \dist(x, \Gamma)^2\big), 
\end{equation*}
where $F(y) = x$. 
The conclusion follows from the definition of $F_*1$ and the fact $F_*I = I$. \proofend

\section{Optimality of the main results}\label{sect-resonant}

In this section, we show that the system is resonant  if the requirements on $A$ and $\Sigma$ mentioned in Theorems~\ref{thm0}, \ref{thm1}, \ref{thm2} are not fulfilled. This implies the optimality of our results. More precisely, we have

\begin{proposition}\label{pro-resonant} Assume that there exists a reflection $F: U \setminus \bar D \to D_\tau$ for some smooth open set $U \subset \subset D$ and some $\tau >0$ such that 
\begin{equation*}
(A, \Sigma) = (F_*A, F_*\Sigma) \mbox{ in } B(x_0, \hat r_0) \cap D, 
\end{equation*}
for some $x_0 \in \Gamma$ and $\hat r_0>0$.  Let $f \in L^2(\mR^d)$ with $\supp f \subset \subset B_{R_0} \setminus \Gamma$ and assume that $A$ is Lipschitz in $\overline{D \cap B(x_0, \hat r_0)}$. There exists $0 < r_0 < \hat r_0$, independent of $f$,  such that if  there is no solution in $H^1(D \cap B(x_0, r_0))$ to the Cauchy problem: 
\begin{equation*}
\dive(A \nabla w) + k^2 \Sigma w = f \mbox{ in } D \cap B(x_0, r_0) \quad \mbox{ and } \quad w = A \nabla w \cdot \nu = 0 \mbox{ on } \partial D \cap B(x_0, r_0),  
\end{equation*}
then   $\limsup_{\delta \to 0}\| u_{\delta}\|_{L^2(K)} = + \infty$ for some $K  \subset \subset B_{R_0} \setminus  \Gamma$ where $u_\delta \in H^1(\mR^d)$ is the unique solution of \eqref{Main-eq-delta}. 
\end{proposition}

Recall that $B(x, r)$ denotes the ball centered at $x$ and of radius $r$. 

\medskip

\noindent{\bf Proof.} Without loss of generality, one may assume that $x_0 = 0$ and $\hat r_0$ is  small. We  prove Proposition~\ref{pro-resonant} by contradiction. Assume that the conclusion is not true. Then even for small $r_0$,  there exists $f$ with $\supp f \cap B_{R_0} \setminus \Gamma$ such that there is no solution in $H^1(D \cap B(x_0, r_0))$ to the Cauchy problem: 
\begin{equation*}
\dive(A \nabla w) + k^2 \Sigma w = f \mbox{ in } D \cap B(x_0, r_0) \quad \mbox{ and } \quad w = A \nabla w \cdot \nu = 0 \mbox{ on } \partial D \cap B(x_0, r_0),  
\end{equation*} 
and 
\begin{equation*}
\limsup_{\delta \to 0} \|u_{\delta} \|_{L^2(K)} < + \infty \mbox{ for all } K \subset \subset B_{R_0} \setminus \Gamma. 
\end{equation*}
Using Lemma~\ref{T0-lem2}, we have 
\begin{equation}\label{hihi}
\| u_{\delta}\|_{H^1(B_{R_0})} \le C \delta^{-1/2}
\end{equation}
since $\supp f \subset \subset B_{R_0} \setminus \Gamma$. 
Set $v_\delta = u_{\delta} \circ F^{-1}$ in $D \cap B(x_0, \hat r_0)$ and  define $w_{\delta} = v_{\delta} - u_{\delta}$ in $D \cap B(x_0, \hat r_0)$. By Lemma~\ref{lem-TO}, we have 
\begin{equation*}
\dive(A \nabla v_{\delta}) + k^2 \Sigma v_{\delta} = -  i \delta_n F_*1 v_{\delta_n} \mbox{ in } D \cap B(x_0, \hat r_0). 
\end{equation*}
Since 
\begin{equation*}
\dive(A \nabla u_{\delta}) + k^2 \Sigma u_{\delta} = k^2 (1 - s_\delta^{-1} s_0 ) \Sigma u_{\delta} -  i \delta s_\delta^{-1}u_{\delta} + s_\delta^{-1} f \mbox{ in } D \cap B(x_0, \hat r_0), 
\end{equation*}
it follows that 
\begin{equation*}
\dive(A \nabla w_{\delta}) + k^2 \Sigma w_{\delta} = g_\delta  \mbox{ in } D \cap B(x_0, \hat r_0), 
\end{equation*}
where 
\begin{equation*}
g_\delta =  f  - i \delta_n F_*1 v_{\delta_n} - k^2 (1 - s_\delta^{-1} s_0 ) \Sigma u_{\delta} +  i \delta s_\delta^{-1}u_{\delta} - (s_\delta^{-1} + 1)  f  \mbox{ in } D \cap B(x_0, \hat r_0). 
\end{equation*}
By Lemma~\ref{lem-TO}, we also have 
\begin{equation*}
w_\delta = 0 \quad \mbox{ and } \quad  A \nabla w_\delta \cdot \nu = i \delta \nabla u_\delta \big|_{D} \cdot \nu \mbox{ on } \partial D \cap  B(x_0, \hat r_0). 
\end{equation*}
Using a local chart and applying  Lemma~\ref{lem-TS} below, we have 
\begin{equation*}
\limsup_{\delta \to 0} \delta^{1/2} \| w_\delta\|_{H^1(D \cap B(x_0, \hat r_0))} = + \infty. 
\end{equation*}
This contradicts \eqref{hihi}. The proof is complete. 
\proofend 

\medskip 
The following lemma is used in the proof of Proposition~\ref{pro-resonant}. 

\begin{lemma}\label{lem-TS} Let  $R > 0$, $a$ be a Lipschitz symmetric uniformly elliptic matrix-valued function and $\sigma$ be a real bounded function  defined  in $B_{R} \cap \mR^d_+$, and let 
$g \in L^2(B_{R})$. Assume that  $W_\delta  \in H^1(B_{R} \cap \mR^d_+)$ $(0 < \delta < 1)$ satisfies
\begin{equation*}
\dive(a \nabla W_\delta) + \sigma W_\delta= g_\delta \mbox{ in } B_{R} \cap \mR^d_+,
\end{equation*}
\begin{equation*}
W_\delta  = 0 \mbox{ on }  B_{R} \cap \mR^d_0, \quad \mbox{ and } \quad   a \nabla W_\delta \cdot \nu   = h_\delta \mbox{ on } B_{R} \cap \mR^d_0, 
\end{equation*}
for some $h_\delta \in H^{-1/2}(B_R \cap \mR^d_0)$ such that 
\begin{equation}\label{ghd}
\|g_\delta - g\|_{L^2(B_R \cap \mR^d_+)} + \| h_\delta \|_{H^{-1/2}(B_R \cap \mR^d_0)} \le c \delta^{1/2},
\end{equation}
for some $c > 0$.  There exists a constant $0 < r < R$ depending only on $R$, and the ellipticity and the Lipschitz  constants of $a$, but  independent of $\delta$, $c$, $g_\delta$, $g$, $h_\delta$, and $\sigma$,  such that if there is {\bf no}  $W \in H^1(B_{r} \cap \mR^d_+)$ with the properties 
\begin{equation}\label{pro-BU-1}
\dive (a \nabla W) = g \mbox{ in } B_{R}  \cap \mR^d_+, \quad W = 0 \mbox{ on }  B_{R}  \cap \mR^d_0, \quad \mbox{and} \quad a \nabla W \cdot \nu = 0  \mbox{ on }  B_{R}  \cap \mR^d_0, 
\end{equation}
then 
\begin{equation}\label{blow-up-1}
\limsup_{\delta \to 0} \delta^{1/2} \| W_\delta\|_{H^1(B_{R}  \cap \mR^d_+)}  = + \infty. 
\end{equation}
\end{lemma}

Here and in what follows, we denote $\mR^d_+ = \mR^d_{e_d,+}$ and  $\mR^d_0 = \mR^d_{e_d,0}$ with $e_d = (0, \cdots, 0, 1) \in \mR^d$. 

\medskip 
\noindent{\bf Proof.} For notational ease, $W_{2^{-n}}$, $g_{2^{-n}}$, and $h_{2^{-n}}$ are denoted  by  $W_n$, $g_n$, and $h_n$ respectively. We have
\begin{equation*}
\dive(a \nabla W_n) = g_n \mbox{ in } B_{R}  \cap \mR^d_+, 
\end{equation*}
\begin{equation*}
W_n  = 0 \mbox{ on } B_{R}  \cap \mR^d_0, \quad  a \nabla W_n \cdot \nu =  h_n  \mbox{ on }B_{R}  \cap \mR^d_0.
\end{equation*}
We prove by contradiction that
\begin{equation}\label{claim-IE}
\limsup_{n \to +\infty} 2^{-n/2}  \| W_n\|_{H^1(B_{R}  \cap \mR^d_+)}   = + \infty. 
\end{equation} 
Assume that 
\begin{equation}\label{claim-IE-contradiction}
m: = \sup_{n}  2^{-n/2}   \| W_n\|_{H^1(B_{R}  \cap \mR^d_+)}   <  + \infty. 
\end{equation} 
Set
\begin{equation*}
w_{n} =\left\{ \begin{array}{cl} W_{n+1} - W_{n} - \bw_n  & \mbox{ in } B_{R}  \cap \mR^d_+\\[6pt]
- \bw_n &  \mbox{ in } B_{R}  \cap \mR^d_-, 
\end{array}\right.
\end{equation*}
where $\bw_n \in H^1(B_R)$ is the unique solution of 
\begin{equation*}
\dive(a \nabla \bw_n) + \sigma \bw_n = (g_{n+1}  - g_n) 1_{B_R \cap \mR^d_+}   \mbox{ in } B_{R} \setminus \mR^d_0, 
\end{equation*}
\begin{equation*}
[a \nabla \bw_n \cdot \nu] = h_{n+1} - h_n \mbox{ on } B_R \cap \mR^d_0, \quad \mbox{ and } \quad a \nabla \bw_n \cdot \nu - i \bw_n = 0 \mbox{ on } \partial B_R. 
\end{equation*}
Here we extend $a$ and $\sigma$ in $B_{R}$ by setting $a(x', x_d) = a(x', -x_d)$ and $\sigma(x', x_d) = 0$ for $(x', x_d) \in (\mR^{d-1} \times \mR_-) \cap B_R$;  though we still denote these extensions by $a$ and $\sigma$. We also denote $1_\Omega$ the characteristic function of a subset $\Omega$ of $\mR^d$. 
We derive from  \eqref{ghd} and \eqref{claim-IE-contradiction} that 
\begin{equation}\label{pro-bwn}
\|\bw_n \|_{H^1(B_R)} \le C m 2^{-n/2}. 
\end{equation}
In this proof, $C$ denotes a constant independent of $n$. From the definition of $w_n$, we have 
\begin{equation*}
\dive(a \nabla w_{n}) + \sigma w_n = 0 \mbox{ in } B_{R}. 
\end{equation*}
From \eqref{claim-IE-contradiction} and \eqref{pro-bwn}, we derive that 
\begin{equation}\label{haha}
\| w_n\|_{H^1(B_R)} \le C m 2^{n/2} \quad \mbox{ and } \quad \|w_n\|_{H^{1}(\partial B_R \cap \mR^d_-)} \le C m 2^{-n/2}. 
\end{equation}
Set $S = (0, \cdots, 0,   -R/4) \in \mR^d$. 
By \cite[Theorem 2]{MinhLoc2} (a three sphere inequality), there exists $r_0 \in (R/4, R/3)$, depending only on $R$ and the Lipschitz and elliptic constants of $a$ such that 
\begin{equation*}
\| w_{n}(\cdot - S) \|_{\bH(\partial B_{r_0}) } \le C  \| w_n (\cdot - S) \|_{\bH(\partial B_{R/4})}^{2/3} \| w_{n} (\cdot - S)\|_{\bH(\partial B_{R/3})}^{1/3}, 
\end{equation*}
where 
\begin{equation*}
\|\varphi \|_{\bH(\partial B_r)}: = \|\varphi \|_{H^{1/2}(\partial B_r)} + \|a \nabla \varphi \cdot \nu \|_{H^{-1/2}(\partial B_r)}. 
\end{equation*}
This implies, by \eqref{haha},  
\begin{equation*}
\| w_{n}(\cdot - S) \|_{\bH(\partial B_{r_0}) } \le  C m 2^{-n/6}.   
\end{equation*}
By Lemma~\ref{T0-lem1}, we obtain 
\begin{equation}\label{ineq-TS}
\| w_{n}(\cdot - S) \|_{H^1(B_{r_0}) } \le C m 2^{-n/6}. 
\end{equation}
Since $\bw_n$ converges in $H^1(B_R)$ by \eqref{pro-bwn}, it follows that $(W_n)$ converges  in $H^1(B_{r} \cap \mR^d_+)$ with $r: = r_0 - R/4$. Let $W$ be the  limit of $W_n$ in $H^1(B_{r} \cap \mR^d_+)$.  Then  
\begin{equation*}
\dive (a \nabla W) = g \mbox{ in } B_{r} \cap \mR^d_+, \quad  W = 0 \mbox{ on }B_{r} \cap  \mR^d_0, \quad  a \nabla W \cdot \eta =  0 \mbox{ on } B_{r} \cap \mR^d_0. 
\end{equation*}
This contradicts the non-existence of $W$.  Hence \eqref{claim-IE} holds. The proof is complete. \proofend

\begin{remark} \fontfamily{m} \selectfont Lemma~\ref{lem-TS} is inspired by \cite[Lemma 2.4]{Ng-CALR}. The proof also has roots from there. The fact that $r$ does not depend on $\sigma$ is somehow surprising. This is based on a new three inequality in \cite[Theorem 2]{MinhLoc2}. Proposition~\ref{pro-resonant} is in the same spirit of the results  in \cite{Ng-CALR} and \cite{KettunenLassas} and extends the results obtained there. 
\end{remark}

\medskip
\noindent {\bf Acknowledgment:} The author thanks Boris Buffoni for interesting discussions on the work of Agmon, Douglis, and Nirenberg in \cite{ADNII}.

\providecommand{\bysame}{\leavevmode\hbox to3em{\hrulefill}\thinspace}
\providecommand{\MR}{\relax\ifhmode\unskip\space\fi MR }
\providecommand{\MRhref}[2]{%
  \href{http://www.ams.org/mathscinet-getitem?mr=#1}{#2}
}
\providecommand{\href}[2]{#2}



\end{document}